\documentclass[10pt,a4paper]{article}
\usepackage{bbm}

\usepackage[leqno]{amsmath}
\usepackage{amsfonts}
\usepackage{graphicx}
\usepackage{amsmath}
\usepackage{amssymb}
\usepackage{latexsym}
\usepackage{amsmath, amsfonts,amssymb, amsthm, euscript,makeidx,color,mathrsfs}
\usepackage{enumerate}

\usepackage[colorlinks,linkcolor=blue,anchorcolor=green,citecolor=red]{hyperref}%when print, use the next package
\def\sqr#1#2{{\vcenter{\vbox{\hrule height.#2pt
              \hbox{\vrule width.#2pt height#1pt \kern#1pt \vrule width.#2pt}
              \hrule height.#2pt}}}}
\def\endpf{\signed {$\sqr69$}}
\def\5n{\negthinspace \negthinspace \negthinspace \negthinspace \negthinspace }
\def\4n{\negthinspace \negthinspace \negthinspace \negthinspace }
\def\3n{\negthinspace \negthinspace \negthinspace }
\def\2n{\negthinspace \negthinspace }
\def\1n{\negthinspace }

\def\dbQ{\mathbb{Q}}
\def\dbR{\mathbb{R}}
\def\dbS{\mathbb{S}}

\def\dbZ{\mathbb{Z}}

\def\={\buildrel \triangle \over =}

\def\ds{\displaystyle}

\def\ns{\noalign{\ss}}
\def\a{\alpha}

\def\d{\delta}
\def\e{\varepsilon}

\def\l{\lambda}
\def\m{\mu}
\def\n{\nu}

\def\t{\tau}
\def\f{\varphi}

\def\ve{\varepsilon}

\def\Th{\Theta}
\def\L{\Lambda}

\def\O{\Omega}
\def\cH{{\cal H}}

\def\cL{{\cal L}}

\def\cU{{\cal U}}
\def\cV{{\cal V}}

\def\be{{\bar e}}

\def\bl{{\bar l}}

\def\br{{\bar r}}
\def\bs{{\bar s}}
\def\bu{{\bar u}}

\def\by{{\bar y}}

\def\hb{{\hat b}}

\def\hy{{\hat y}}

\def\tix{{\tilde x}}
\def\tiy{{\tilde y}}
\def\tiz{{\tilde z}}

\def\timu{{\tilde \m}}

\def\ss{\smallskip}
\def\ms{\medskip}
\def\bs{\bigskip}
\def\q{\quad}
\def\qq{\qquad}
\def\hb{\hbox}

\def\da{\mathop{\downarrow}}
\def\Ra{\mathop{\Rightarrow}}

\def\lan{\langle}
\def\ran{\rangle}

\def\rf{\eqref}
\def\h{\widehat}
\def\wt{\widetilde}

\def\cd{\cdot}
\def\cds{\cdots}

\def\ae{\hbox{\rm a.e.}}

\def\les{\leqslant}
\def\ges{\geqslant}

\def\({\Big (}
\def\){\Big )}
\def\[{\Big[}
\def\]{\Big]}
\def\pa{{\partial}}
\def\bde{\begin{definition}\label}
\def\ede{\end{definition}}
\def\be{\begin{equation}}
\def\bel{\begin{equation}\label}
\def\ee{\end{equation}}
\def\bt{\begin{theorem}\label}
\def\et{\end{theorem}}
\def\bc{\begin{corollary}\label}
\def\ec{\end{corollary}}
\def\bl{\begin{lemma}\label}
\def\el{\end{lemma}}
\def\bp{\begin{proposition}\label}
\def\ep{\end{proposition}}
\def\bas{\begin{assumption}\label}
\def\eas{\end{assumption}}
\def\br{\begin{remark}\label}
\def\er{\end{remark}}
\def\bex{\begin{example}\label}
\def\ex{\end{example}}
\def\ba{\begin{array}}
\def\ea{\end{array}}
\def\ed{\end{document}}
\def\square#1{\vbox{\hrule\hbox{\vrule height#1%
     \kern#1\vrule}\hrule}}
\def\rectangle#1#2{\vbox{\hrule\hbox{\vrule height#1%
     \kern#2\vrule}\hrule}}

\font\tenbb=msbm10 \font\sevenbb=msbm7 \font\fivebb=msbm5

\newfam\bbfam
\scriptscriptfont\bbfam=\fivebb \textfont\bbfam=\tenbb
\scriptfont\bbfam=\sevenbb

\newtheorem{theorem}{\hskip 1.3em Theorem}[section]
\newtheorem{definition}[theorem]{\hskip 1.3em Definition}
\newtheorem{proposition}[theorem]{\hskip 1.3em Proposition}
\newtheorem{corollary}[theorem]{\hskip 1.3em Corollary}
\newtheorem{lemma}[theorem]{\hskip 1.3em Lemma}
\newtheorem{remark}[theorem]{\hskip 1.3em Remark}
\newtheorem{example}[theorem]{\hskip 1.3em Example}

\newtheorem{assumption}[theorem]{\hskip 1.3em Assumption}
\makeatletter
   
   \@addtoreset{equation}{section}
\makeatother

\def\endpf{\hfill$\Box$\vspace{0.4cm}\\ }

\def\eqon{ \, {\rm on } \,\, }
\def\eqin{ \, {\rm in } \,\, }

\def\Div{ \, {\nabla\cd\1n} \,\, }

\begin{document}

\title{\bf Second-Order Necessary Conditions for \\
 Optimal Control of Semilinear
 Elliptic Equations \\ with Leading Term Containing Controls
\thanks{The first author was supported in part by NSFC grant 11371104,
 the second author was supported in part by NSF grant DMS-1406776.}
\author{\rm Hongwei Lou\footnote{School of Mathematical Sciences, and LMNS, Fudan
University, Shanghai 200433, China (Email:
\texttt{hwlou@fudan.edu.cn})}~~~and~~~Jiongmin
Yong\footnote{Department of Mathematics, University of Central
Florida, Orlando, FL 32816, USA (Email: \texttt{jiongmin.yong@ucf.edu})}}}

\date{}
\maketitle
\begin{quote}
%\footnotesize

\vskip-0.8cm

%\centerline{({\it Dedicated to Professor Eduardo Casas at his 60th birthdate})}

\bs

{\bf Abstract.} An optimal control problem for a semilinear elliptic equation of divergence
form is considered. Both the leading term and the semilinear term of the state equation
contain the control. The well-known Pontryagin type maximum principle for the optimal controls is the {\it first}-order necessary condition. When such a first-order necessary condition is singular in some sense, certain type of the {\it second}-order necessary condition will come in naturally. The aim of this paper is to explore such kind of conditions for our optimal control problem.

\ms


\ms


\end{quote}

\normalsize

\newtheorem{Definition}{Definition}[section]
\newtheorem{Theorem}[Definition]{Theorem}
\newtheorem{Lemma}[Definition]{Lemma}
\newtheorem{Corollary}[Definition]{Corollary}
\newtheorem{Proposition}[Definition]{Proposition}
\newtheorem{Remark}{Remark}[section]
\newtheorem{Example}{Example}
%\newfont{\Bbb}{msbm10 scaled\magstephalf}
%\newfont{\frak}{eufm10 scaled\magstephalf}
%\newfont{\sfr}{msbm7 scaled\magstephalf}

%\def\theequation{1.\arabic{equation}}

\section{Introduction.}

Let $\O\subseteq\dbR^n$ ($n\ges2$) be a bounded domain with a smooth boundary $\pa\O$.
Consider the following controlled elliptic partial differential equation (PDE, for short):
\bel{state}\left\{\2n\ba{ll}
\ns\ds-\Div\big(A(x,u(x))\nabla
y(x)\big)=f(x,y(x),u(x)),\qq\hb{in }~\O,\\
\ns\ds y(x)=0,\qq\qq\qq\qq\qq\hb{on }~\pa\O,\ea\right.\ee
where $A:\O\times U\to\dbS^n_+$ and $f:\O\times\dbR\times U\to\dbR$, with $\dbS^n_+$ being the set of all $(n\times n)$ positive definite matrices, $U$ being a separable (nonempty) metric space. In the above, $u(\cd)$ is the {\it control} which belongs to the set $\cU$ of all {\it admissible controls} defined by the following:
$$\cU\equiv\big\{u:\O\to U\bigm|u(\cd)\hb{ is measurable}\,\big\}.$$
Under some mild conditions, for any $u(\cd)\in\cU$, \rf{state} admits a unique weak solution $y(\cd)\equiv y(\cd\,;u(\cd))$ which is called the {\it state} (corresponding to the control $u(\cd)$). The performance of the control $u(\cd)$ is measured by the following cost functional
\bel{J(u)}J(u(\cd))=\int_\O f^0(x,y(x),u(x))dx\equiv\int_\O f^0(x,y(x;u(\cd)),u(x))dx\ee
for some given map $f^0:\O\times\dbR\times U\to\dbR$. Our optimal
control problem can be stated as follows.

\ms

{\bf Problem (C)}. Find a $\bu(\cd)\in\cU$ such that
\bel{J=inf}J(\bu(\cd))=\inf_{u(\cd)\in\cU}J(u(\cd)).\ee
Any $\bar u(\cd)\in\cU$ satisfying \rf{J=inf} is called an {\it
optimal control}, and the corresponding $\by(\cd)\equiv
y(\cd\,;\bar u(\cd))$ is called an {\it optimal state}. The pair
$(\bar y(\cd),\bar u(\cd))$ is called an {\it optimal pair}.

\ms

Let us make some rough observations. Suppose $(\bar y(\cd),\bar u(\cd))$ is an optimal pair of Problem (C). For any given $u(\cd)\in\cU$, let $u^\d(\cd)\in\cU$ be a suitable perturbation of $\bar u(\cd)$ determined by $u(\cd)$ with a parameter $\d>0$ (for examples, a convex type perturbation, or a spike type variation), so that $\wt\rho(u^\d(\cd),\bar u(\cd))=O(\d)$ with $\wt\rho(\cd\,,\cd)$ being a suitable metric on the set $\cU$, and the following holds:
\bel{J=J+dJ^1}J(u^\d(\cd))=J(\bar u(\cd))+\d J^1(\bar u(\cd),u(\cd))+o(\d),\qq\hb{as }\d\to0.\ee
Here, $J^1(u(\cd),u(\cd))$ is some functional of $(\bar u(\cd),u(\cd))$. The above can be called the first-order Taylor expansion of $J(\cd)$ at $\bar u(\cd)$, and $J^1(\bar u(\cd),u(\cd))$ can be regarded as the ``directional derivative'' of $J(\cd)$ at $\bar u(\cd)$ in the ``direction'' $u(\cd)$. Hence, the minimality of $\bar u(\cd)$ implies
\bel{J1>0}J^1(\bar u(\cd),u(\cd))\ges0,\qq\forall u(\cd)\in\cU.\ee
This is called the {\it first-order necessary condition} for $\bar u(\cd)$, which is essentially the Pontryagin's maximum principle for our Problem (C). Now, suppose that there is a set $\cU_0\subseteq\cU$, which is different from the singleton $\{\bar u(\cd)\}$, such that the following holds:
\bel{J^1=0}J^1(\bar u(\cd),u(\cd))=0,\qq\forall u(\cd)\in\cU_0.\ee
Then $\bar u(\cd)$ is said to be {\it singular} on the set $\cU_0$. For convenience, we call $\cU_0$ a {\it singular set} of $\bar u(\cd)$. Let
$$\cU_0(\bar u(\cd))=\Big\{u(\cd)\in\cU\bigm|J^1(\bar u(\cd),u(\cd))=0\Big\},$$
which is called the {\it maximum singular set} of $\bar u(\cd)$. When $\cU_0(\bar u(\cd))=\cU$,
we say that $\bar u(\cd)$ is {\it fully singular} (or simply {\it singular}); When $\cU_0(\bar u(\cd))=\{\bar u(\cd)\}$, we say that $\bar u(\cd)$ is {\it nonsingular}; And, more interestingly, when $\{\bar u(\cd)\}\ne\cU_0(\bar u(\cd))\ne\cU$, we say that $\bar u(\cd)$ is {\it partially singular}. The notion of singular control was introduced by Gabasov--Kirillova in \cite{Gabasov-Kirillova 1972}, where our partial singularity was called ``the singularity in the sense of Pontryagin's maximum principle'', and our full singularity was called ``the singularity in the classical sense''. We prefer to use our shorter names. Now, suppose $\bar u(\cd)$ is partially singular. Then one should expect that the following (comparing with \rf{J=J+dJ^1})
\bel{J=J+d^2J^2}J(u^\d(\cd))=J(\bar u(\cd))+\d^2J^2(\bar u(\cd),u(\cd))+o(\d^2),\qq\forall u(\cd)\in\cU_0(\bar u(\cd)),\ee
for some functional $J^2(\bar u(\cd),u(\cd))$ of $(\bar u(\cd),u(\cd))$. The above can be called the second-order Taylor expansion of $J(\cd)$ at $\bar u(\cd)$ in the direction of $u(\cd)\in\cU_0(\bar u(\cd))$, and $J^2(\bar u(\cd),u(\cd))$ can be regarded as the ``second order directional derivative'' at $\bar u(\cd)$ in the ``direction'' of $u(\cd)\in\cU_0(\bar u(\cd))$. Then the minimality of $\bar u(\cd)$ leads to the following:
\bel{J^2>0}J^2(\bar u(\cd),u(\cd))\ges0,\qq\forall u(\cd)\in\cU_0(\bar u(\cd)).\ee
This is referred to as the {\it second-order necessary condition} of $\bar u(\cd)$. We emphasize that the above holds only for all $u(\cd)\in\cU_0(\bar u(\cd))$, the maximum singular set of $\bar u(\cd)$. To get some more feeling, let us look at the following simple example, consisting of three situations.

\bex{1.1} \rm Let $U=[-1,1]\times[-1,1]$.

\ms

(i) Let $J(u)=u^2=u_1^2+u_2^2$ with $u=(u_1,u_2)\in U$. Then $u\mapsto J(u)$ is differentiable and the minimum is attained at $\bar u=(0,0)$, an interior point of $U$, with $J(\bar u)=0$. Therefore, for any $u=(u_1,u_2)\in U$ and $\d\in(0,1)$, we have $u^\d=\bar u+\d(u-\bar u)=\d u\in U$, and
$$J(u^\d)=\d^2u^2\equiv J(\bar u)+\d J^1(\bar u,u)+\d^2J^2(\bar u,u).$$
Consequently,
$$J^1(\bar u,u)=0,\qq\forall u\in U.$$
This means that the maximum singular set $U_0(\bar u)$ of $\bar u$ coincides with $U$, and $\bar u$ is fully singular. Hence,
$$J^2(\bar u,u)=u^2\ges0,\qq\forall u\in U_0(\bar u)=U,$$
which is the classical second-order necessary condition for $\bar u$.

\ms

(ii) Let $J(u)=u_1^3+u_2^3$ with $u=(u_1,u_2)\in U$. The minimum is attained at $\bar u=(-1,-1)$. Then, for any $u=(u_1,u_2)\in U$ and any $\d\in(0,1)>0$, one has
$$u^\d=\bar u+\d(u-\bar u)=(-1,-1)+\d(u_1+1,u_2+1)\in U,$$
and
$$\ba{ll}
\ns\ds J(u^\d)=\big[-1+\d(u_1+1)\big]^3+\big[-1+\d(u_2+1)\big]^3\\
\ns\ds\qq~=-2+3\d(u_1+u_2+2)-3\d^2[(u_1+1)^2+(u_2+1)^2]+\d^3[(u_1+1)^3+(u_2+1)^3]\\
\ns\ds\qq~=J(\bar u)+\d J^1(\bar u,u)+\d^2J^2(\bar u,u)+\d^3J^3(\bar u,u).\ea$$
Thus,
$$J^1(\bar u,u)=u_1+u_2+2\ne0,\qq\forall u=(u_1,u_2)\in U\setminus\{\bar u\}.$$
This means that $\bar u$ is nonsingular. In this case, there is no second-order necessary condition for $\bar u$.

\ms

(iii) Let $J(u)=u_1^3+u_2^2$ with $u=(u_1,u_2)\in U$. The minimum is attained at $\bar u=(-1,0)$. For any $u=(u_1,u_2)\in U$, let the perturbation $u^\d$ be defined by the following:
$$u^\d=\bar u+\d(u-\bar u)=(-1,0)+\d(u_1+1,u_2)=(-1+\d(u_1+1),\d u_2)\in U.$$
Then we have
\bel{1.9}\ba{ll}
\ns\ds J(u^\d)=[-1+\d(u_1+1)]^3+\d^2u_2^2=J(\bar u)+3\d(u_1+1)
+\d^2[-3(u_1+1)^2+u_2^2]+\d^3(u_1+1)^3\\
\ns\ds\qq\q\equiv J(\bar u)+\d J^1(\bar u,u)+\d^2J^2(\bar u,u)+\d^3J^3(\bar u,u).\ea\ee
Hence, the first-order necessary condition is
$$J^1(\bar u,u)\equiv3(u_1+1)\ges0,\qq\forall u=(u_1,u_2)\in U,$$
and $\bar u$ is partially singular with $U_0(\bar u)=\{(-1,u_2)\bigm|u_2\in[-1,1]\}$. The second-order necessary condition is
$$J^2(\bar u,u)\equiv u_2^2\ges0,\qq\forall u\in U_0(\bar u).$$
However, we do not have (see \rf{1.9})
$$J^2(\bar u,u)\equiv-3(u_1+1)^2+u_2^2\ges0,\qq\forall u\in U.$$

\ex

The above example shows that in general, fully singular, partially singular, and nonsingular all can happen for a minimum of a function. Of course, the above is for scalar functions, and it is expected that the case of optimal control problems should be much more complicated. In Lou \cite{Lou 2010}, second-order necessary and sufficient conditions for partially singular optimal controls of ordinary differential equations were established with general control domain. On the other hand, for convex control domains (mainly interval type), and mainly for fully singular cases, second-order necessary/sufficient optimality conditions have been studied for PDEs by many authors. We mention just a few of them here: Casas--Tr\"oltzsch \cite{Casas-Troltzsch 1999, Casas-Troltzsch 2002, Casas-Troltzsch 2009}, Casas--Tr\"oltzsch--Unger \cite{Casas-Troltzsch-Unger 2000}, Raymond--Tr\"oltzsch \cite{Raymond-Troltzsch 2000}, Mittelmann \cite{Mittelmann 2001}, Casas--Mateos \cite{Casas-Mateos 2002}, R\"osch--Tr\"oltzsch \cite{Rosch-Troltzsch 2003, Rosch-Troltzsch 2006}, Wang--He \cite{Wang-He 2006}, Casas--Los Reye--Tr\"oltzsch \cite{Casas-Los Reye-Troltzsch 2008}, and Bonnans--Hermant \cite{Bonnans-Hermant 2009a, Bonnans-Hermant 2009b}. For some earlier works on ODEs, see Kelly \cite{Kelly 1964}, Kopp--Moyer \cite{Kopp-Moyer 1965}, Gabasov--Kirillova \cite{Gabasov-Kirillova 1972}, Krener \cite{Krener 1977}, and Knobloch \cite{Knobloch 1981}.

\ms

For problems of elliptic PDEs with control appearing in the leading term, Casas \cite{Casas 1992} studied the first-order necessary conditions for the case $A(x,u)=uI$ with quadratic cost functional and with the control being Lipschtz continuous. General case were treated by Lou--Yong in \cite{Lou-Yong 2009}, and analogous results for parabolic and hyperboliccases were given by Lou in \cite{Lou 2011} and Li-Lou in \cite{Li-Lou 2012}. If the leading term of the state equation \rf{state} does not contain controls, i.e., $A(x,u)\equiv A(x)$, then one can establish the second-order necessary conditions for partially singular optimal controls following similar arguments of \cite{Lou 2010}. However, if the leading term of the equation contains the control, we will see that it is much more complicated, even in defining the partial singularity of the optimal control. It turns out that the construction of a proper family of perturbations is much more difficult than the case without having control in the leading term, in order to have the first-order term disappeared in the Taylor type expansion. The difficult will be overcome by introducing the notion of weak singularity which involves a proper vector field. Consequently, the results obtained will have some big difference comparing with those for the problems without having the control in the leading term.

\ms

The rest of the paper is organized as follows. In Section 2, we will introduce the notions of singularity and weak singularity of the optimal controls. The main result of the paper will be stated, together with a couple of corollaries. Section 3 will be devoted to a review of the proof for the first-order necessary condition for Problem (C), which will inspire the second-order necessary condition. Section 4 is devoted to a proof of a result crucial for the proof our main result. A proof of the second-order necessary condition will be presented in Section 5.

\section{The Main Result.}

For any differentiable function $\f:\O\to\dbR$, its gradient is denoted by $\nabla\f=({\pa\f\over\pa x_1},\ldots,{\pa\f\over\pa x_n})^\top:\O\to\dbR^n(\equiv\dbR^{n\times1})$; For any differentiable vector-valued function $f=(f^1,f^2,\ldots,f^n)^\top:\O\to\dbR^n$, its Jacobean matrix $f_x$ is denoted by
$$f_x=\begin{pmatrix}{\pa f^1\over\pa x_1}&{\pa
f^2\over\pa x_1}&\ldots&{\pa f^n\over\pa x_1}\\ \\
{\pa f^1\over\pa x_2}&{\pa f^2\over\pa x_2}&\ldots&{\pa f^n\over\pa x_2}\\
\vdots&\vdots&\ddots&\vdots\\
{\pa f^1\over\pa x_n}&{\pa f^2\over\pa x_n}&\ldots&{\pa
f^n\over\pa x_n}\end{pmatrix}\equiv\(\nabla f^1,\nabla f^2,\cds,\nabla f^n\)\equiv\nabla f^\top.$$
Compatible with the above notation, we also will use
$$\nabla^\top\1n F=(\nabla \cd F^1, \nabla\cd F^2,\ldots, \nabla \cd F^n).$$
for $F=(F^1, F^2,\ldots, F^n):\O\to\dbR^{n\times n}$.
\if{
$$\nabla^\top\2n f=\big(f_x\big)^\top\2n=\begin{pmatrix}{\pa f^1\over\pa x_1}&{\pa
f^1\over\pa x_2}&\ldots&{\pa f^1\over\pa x_n}\\ \\
{\pa f^2\over\pa x_1}&{\pa f^2\over\pa x_2}&\ldots&{\pa f^2\over\pa x_n}\\
\vdots&\vdots&\ddots&\vdots\\
{\pa f^n\over\pa x_1}&{\pa f^n\over\pa x_2}&\ldots&{\pa
f^n\over\pa x_n}\end{pmatrix}\equiv\begin{pmatrix}(\nabla f^1)^\top\\ (\nabla f^2)^\top\\
\vdots\\ (\nabla f^n)^\top\end{pmatrix}.$$
}\fi
A function $g:\dbR^n\to\dbR$ is said to be $[0,1]^n$--periodic if it admits a period $1$ in every coordinate direction $x_i$, $i=1,2,\cds,n$. Denote $W^{1,2}_\#([0,1]^n;\dbR^n)$ the space of all $[0,1]^n$--periodic vector-valued functions in $W^{1,2}_{loc}(\dbR^n;\dbR^n)$ and $W^{1,2}_\#([0,1]^n;\dbR^n)/\dbR^n$ the corresponding quotient space.

\ms

Next, let us introduce the following assumptions.

\ms

{\bf(S1)} Set $\O$ is a bounded domain in $\dbR^n$ ($n\ges2$) with a smooth boundary $\pa\O$, and metric space $(U,\rho)$ is separable.

\ms

{\bf(S2)} Function $A:\O\times U\to\dbS^n_+$, (recall that $\dbS^n_+$ is the set of all $(n\times n)$ (symmetric) positive definite matrices), for which $x\mapsto A(x,v)$ is measurable, and  $v\mapsto A(x,v)$ is continuous. Further, there exist constants $\L\ges\l>0$ such that
\bel{l<A<L}\l|\xi|^2\les\lan A(x,v)\xi,\xi\ran\les\L|\xi|^2,\qq\forall
\xi\in\dbR^n,\,\ae\,x\in\O,\,v\in U.\ee

\ms

{\bf(S3)} Function $f:\O\times\dbR\times U\to\dbR$ has the following properties: $x\mapsto f(x,y,v)$ is measurable, $(y,v)\mapsto f(x,y,v)$ is continuous for almost all $x\in\O$, and $y\mapsto f(x,y,v)$ continuously differentiable. Moreover,
\bel{f_y<0}f_y(x,y,v)\les0,\qq\ae\,(x,y,v)\in\O\times\dbR\times U\ee
and for any $R>0$, there exists an $M_R>0$ such that
\bel{|f|<M}|f(x,y,v)|+|f_y(x,y,v)|\les M_R,\qq\ae\,(x,v)\in\O\times U,~|y|\les R.\ee

\ms

{\bf(S4)} Function $f^0:\O\times\dbR\times U\to\dbR$ has the following properties: $x\mapsto f^0(x,y,v)$ is measurable, $(y,v)\mapsto f^0(x,y,v)$ is continuous for almost all $x\in\O$, and $y\mapsto f^0(x,y,v)$ is continuously differentiable. Moreover, for any $R>0$, there exists a $K_R>0$ such that
\bel{|f^0|<M}|f^0(x,y,v)|+|f^0_y(x,y,v)|\les K_R,\qq\ae\,(x,v)\in\O\times U, \, |y|\leq R.\ee

\ms

It is standard that under (S1)--(S3), for any $u(\cd)\in\cU$, state equation \rf{state} admits a unique weak solution $y(\cd)=y(\cd\,;u(\cd))\in H_0^1(\O)\cap C(\bar\O)$ and the following estimate holds:
\bel{|y|<K}\|y(\cd)\|_{H^1_0(\O)}+\|y(\cd)\|_{L^\infty(\O)}\les K,\ee
for some constant $K>0$. Therefore, if, in addition, (S4) is also assumed, then the cost functional is well-defined. Consequently, Problem (C) is well-formulated. The following was established in \cite{Lou-Yong 2009}.

\bt{T201}  \sl Let {\rm(S1)--(S4)} hold. Let $(\bar y(\cd),\bar u(\cd))$ be an optimal pair of Problem {\rm(C)}, and $\bar\psi(\cd)$ be the weak solution of the following {\it adjoint equation}:
\bel{adjoint}\left\{\2n\ba{ll}
\ns\ds-\Div\big(A(x,\bu(x))\nabla\bar\psi(x)\big)=f_y(x,\bar y(x),\bar u(x))\,\bar\psi(x)
-f^0_y(x,\bar y(x),\bar u(x)),\qq\eqin  \O,\\ [2mm]
\ds\bar\psi\big|_{\,\pa\O}=0.\ea\right.\ee
Then
\bel{H-H>0}\ba{ll}
\ns\ds H\big(x,\bar y(x),\bar\psi(x),\nabla\bar y(x),\nabla\bar\psi(x),\bar u(x)\big)
-H\big(x,\bar y(x),\bar\psi(x),\nabla\bar y(x),\nabla\bar\psi(x),v\big)\\
\ns\ds\ges\max_{\m\in S^{n-1}}{\lan\big[A(x,\bar u(x))-A(x,v)\big]\nabla\bar y(x),\m\ran\,
\lan\big[A(x,\bar u(x))-A(x,v)\big]\nabla\bar\psi(x),\m\ran\over\lan A(x,v)\m,\m\ran}\\
\ns\ds={1\over 2}\big|A(x,v)^{-{1\over
2}}(A(x,\bar u(x))-A(x,v))\nabla\bar y(x)\big|\,\big|A(x,v)^{-{1\over
2}}(A(x,\bu(x))-A(x,v))\nabla\bar\psi(x)\big|\\
\ns\ds\qq+{1\over 2}\lan A(x,v)^{-{1\over2}}(A(x,\bar u(x))-A(x,v))\nabla\bar y(x),A(x,v)^{-{1\over2}}(A(x,\bar u(x))-A(x,v))\nabla\bar\psi(x)\ran\ges0,\\
\ns\ds\qq\qq\qq\qq\qq\qq\qq\qq\qq\qq\qq\qq\forall v\in U,\q\ae\,x\in\O,\ea\ee
where
\bel{H}\ba{ll}
\ns\ds H(x,y,\psi,\xi,\eta,v)=\psi f(x,y,v)-f^0(x,y,v)-\lan A(x,v)\xi,\eta\ran,\\
\ns\ds\qq\qq\qq\qq\qq\forall (x,y,\psi,\xi,\eta,v)\in\dbR^n\times\dbR\times \dbR\times\dbR^n\times\dbR^n\times U,\ea\ee
which is called the {\it Hamiltonian}, and $S^{n-1}$ is the unit sphere in $\dbR^n$.

\et

The equality in \rf{H-H>0} follows from the following simple fact (see \cite{Lou-Yong 2009}, Lemma 2.3, for a proof).
\bel{max-xi-eta}\max_{\m\in S^{n-1}}\lan\m,\xi\ran\,\lan\m,\eta\ran={|\xi|\,|\eta|+\lan\xi,\eta\ran\over 2},\qq\forall \xi,\eta\in\dbR^n,~n\ges 2,\ee
and the last inequality in \rf{H-H>0} is due to the Cauchy-Schwartz inequality. Therefore, \rf{H-H>0} implies (and might be a little stronger than) the following:
\bel{H=max H} H\big(x,\bar y(x),\bar\psi(x),\nabla\bar y(x),\nabla\bar\psi(x),\bar u(x)\big)=
\max_{v\in U}H\big(x,\bar y(x),\bar\psi(x),\nabla\bar y(x),\nabla\bar\psi(x),v\big),\q\ae\, x\in\O.\ee

\ms

Let
$$\cL=\Big\{\ell:\O\to S^{n-1}\bigm|\ell(\cd)\hb{ is measurable }\Big\}.$$
We now introduce the following definition.

\bde{T204} \rm Let $\bar u(\cd)\in\cU$ be an optimal control of Problem {\rm(C)}.

\ms

(i) Let $(u(\cd),\ell(\cd))\in\cU\times\cL$ satisfy the following:
\bel{weakly singular}\ba{ll}
\ns\ds H\big(x,\bar y(x),\bar\psi(x),\nabla\bar y(x),\nabla\bar\psi(x),\bar u(x)\big)
-H\big(x,\bar y(x),\bar\psi(x),\nabla\bar y(x),\nabla\bar\psi(x),u(x)\big)\\
\ns\ds={\lan\big[A(x,\bar u(x))-A(x,u(x))\big]\nabla\bar y(x),\ell(x)\ran\,\lan
\big[A(x,\bar u(x))-A(x,u(x))\big]\nabla\bar\psi(x),\ell(x)\ran\over
\lan A(x,u(x))\ell(x),\ell(x)\ran},\q\ae\,x\in\O.\ea\ee
Then we say that $\bar u(\cd)$ is {\it weakly singular} at $(u(\cd),\ell(\cd))$.

\ms

(ii) Let $u(\cd)\in\cU$ such that
\bel{singular}H\big(x,\bar y(x),\bar\psi(x),\nabla\bar y(x),\nabla\bar\psi(x),\bar u(x)\big)=
H\big(x,\bar y(x),\bar\psi(x),\nabla\bar y(x),\nabla\bar\psi(x),u(x)\big),\q\ae\,x\in\O,\ee
then we say that $\bar u(\cd)$ is {\it singular} at $u(\cd)$.

\ms

(iii) Denote
$$\ba{ll}
\ns\ds\cV_0(\bar u(\cd))=\Big\{(u(\cd),\ell(\cd))\in\cU\times\cL\bigm|\bar u(\cd)\hb{ is weakly singular at }(u(\cd),\ell(\cd))\Big\},\\
\ns\ds\cU_0(\bar u(\cd))=\Big\{u(\cd)\in\cU\bigm|\bar u(\cd)\hb{ is singular at }u(\cd)\Big\}.\ea$$
If $\cV_0(\bar u(\cd))=\cU\times\cL$, we say that $\bar u(\cd)$ is {\it fully weakly singular}; If
$$\cV_0(\bar u(\cd))\ne\cU\times\cL,\qq\cV_0(\bar u(\cd))\setminus\(\{\bar u(\cd)\}\times\cL\)\ne\emptyset,$$
then we say that $\bar u(\cd)$ is {\it partially weakly singular}; If $\cV_0(\bar u(\cd))\subseteq\{\bar u(\cd)\}\times\cL$ (in this case, the equality actually holds), then we say that $\bar u(\cd)$ is {\it weakly nonsingular}.

\ms

Likewise, we may define $\bu(\cd)$ to be {\it fully singular}, {\it partially singular}, and {\it  nonsingular}, respectively, when $\cU_0(\bu(\cd))=\cU$, $\{\bar u(\cd)\}\ne\cU_0(\bu(\cd))\ne\cU$ and
$\cU_0(\bu(\cd))=\{\bu(\cd)\}$, respectively.

\ede

Let us make some observations on the above notions.

\ms

$\bullet$ If optimal control $\bar u(\cd)$ is weakly singular at $(u(\cd),\ell(\cd))$ (with $u(\cd)\ne\bu(\cd)$), then comparing \rf{weakly singular} with \rf{H-H>0}, we see that for almost all $x\in\O$, $(u(x),\ell(x))$ is a maximum of the map
$$(v,\m)\mapsto{\lan\big[A(x,\bar u(x))-A(x,v)\big]\nabla\bar y(x),\m\ran\,
\lan\big[A(x,\bar u(x))-A(x,v)\big]\nabla\bar\psi(x),\m\ran\over\lan A(x,v)
\m,\m\ran}$$
over $U\times S^{n-1}$. Note that such a maximum point might not be unique, in general.

\ms

$\bullet$ If optimal control $\bu(\cd)$ is singular at $u(\cd)\ne\bu(\cd)$, then by \rf{H-H>0}, we see that
\bel{max=0}\ba{ll}
\ns\ds\max_{\m\in S^{n-1}}{\lan\big[A(x,\bar u(x))-A(x,u(x))\big]\nabla\bar y(x),\m\ran\,
\lan\big[A(x,\bar u(x))-A(x,u(x))\big]\nabla\bar\psi(x),\m\ran\over\lan A(x,u(x))\m,\m\ran}\\
\ns\ds={1\over 2}\big|A(x,u(x))^{-{1\over
2}}[A(x,\bar u(x))-A(x,u(x))]\nabla\bar y(x)\big|\,\big|A(x,u(x))^{-{1\over
2}}[A(x,\bu(x))-A(x,u(x))]\nabla\bar\psi(x)\big|\\
\ns\ds~+\1n{1\over 2}\lan A(x,\1n u(x))^{-{1\over2}}[A(x,\1n\bar u(x))\1n-\1n A(x,\1n u(x))]\nabla\bar y(x),A(x,\1n u(x))^{-{1\over2}}[A(x,\1n\bar u(x))\1n-\1n A(x,\1n u(x))]\nabla\bar\psi(x)\ran\1n=\1n0,\\
\ns\ds\qq\qq\qq\qq\qq\qq\qq\qq\qq\qq\qq\qq\qq\qq\ae\,x\in\O.\ea\ee
Therefore, by the compactness of $S^{n-1}$, together with Filipov's measurable selection lemma (\cite{Li-Yong 1995}), we have some $\ell(\cd)\in\cL$ such that
\bel{(A-A)(A-A)=0}\lan\big[A(x,\bar u(x))-A(x,u(x))\big]\nabla\bar y(x),\ell(x)\ran\,
\lan\big[A(x,\bar u(x))-A(x,u(x))\big]\nabla\bar\psi(x),\ell(x)\ran=0,\q\ae~x\in\O.\ee
Hence, the optimal control $\bar u(\cd)$ is weakly singular at $(u(\cd),\ell(\cd))$ for some  $\ell(\cd)\in\cL$.

\ms

$\bullet$ If optimal control $\bu(\cd)$ is weakly singular at $(u(\cd),\ell(\cd))$ such that \rf{(A-A)(A-A)=0} holds, then $\bu(\cd)$ must be singular at $u(\cd)$. Further, it follows from the last equality in \rf{max=0} that the equality holds in Cauchy-Schwartz inequality. Therefore,
$$[A(x,\bar u(x))-A(x,u(x))]\nabla\bar y(x)\q\hb{and}\q[A(x,\bar u(x))-A(x,u(x))]\nabla\bar\psi(x)$$
must be linearly dependent and have opposite directions. Consequently, \rf{(A-A)(A-A)=0} implies that
\bel{(A-A)=0}\lan\big[A(x,\bar u(x))-A(x,u(x))\big]\nabla\bar y(x),\ell(x)\ran=
\lan\big[A(x,\bar u(x))-A(x,u(x))\big]\nabla\bar\psi(x),\ell(x)\ran=0,\q\ae~x\in\O.\ee
The above can be summarized as follows.

\bp{P203} \sl Let {\rm(S1)--(S4)} hold. Suppose $\bu(\cd)$ is an optimal control of Problem {\rm(C)}. If $\bu(\cd)$ is singular at $u(\cd)\in\cU\setminus\{\bu(\cd)\}$. Then there exists an $\ell(\cd)\in\cL$ such that \rf{(A-A)=0} holds and $\bu(\cd)$ is weakly singular at $(u(\cd),\ell(\cd))$. Conversely, if $\bu(\cd)$ is weakly singular at $(u(\cd),\ell(\cd))\in\cU\times\cL$ such that \rf{(A-A)=0} holds, then $\bu(\cd)$ is singular at $u(\cd)$.

\ep

If $A(x,v)$ is independent of $v\in U$, then the right hand side of \rf{weakly singular} is automatically zero, and \rf{singular} is true. Thus, in such a case, weak singularity is equivalent to singularity, and
$$\cV_0(\bar u(\cd))=\cU_0(\bu(\cd))\times\cL.$$

\ms

To state our main result of the current paper, the second-order necessary condition for optimal control of Problem (C), we need the following further assumption.

\ms

{\bf(S5)} Function $y\mapsto(f(x,y,v),f^0(x,y,v))$ is twice continuously differentiable.
Moreover, for any $R>0$, there exists a $K_R>0$ such that
\bel{|f_yy|<K}|f_{yy}(x,y,v)|+|f^0_{yy}(x,y,u)|\les K_R,\qq\forall v\in U,\q|y|\les R,\q\ae\,x\in\O.\ee

\ms

We point out that, unlike most of the literature on PDE controls that we cited, no differentiability condition is assumed for the map $u\mapsto(f(x,y,u),f^0(x,y,u))$. Actually, our $U$ is just a metric space which does not have a linear structure, in general. In particular, no convexity condition is assumed for $U$. Now, we state our main result of this paper.

\bt{T204} \sl Let {\rm(S1)--(S5)} hold and $(\bar y(\cd),\bar u(\cd))$ be an optimal pair of Problem {\rm(C)}. Let $\bar u(\cd)$ be partially weakly singular and $(u(\cd),\ell(\cd))\in\cV_0(\bar u(\cd))$ with $u(\cd)\ne\bar u(\cd)$. Then the following holds:
\bel{second}\ba{ll}
\ns\ds\int_\O\2n\Big\{\1n
\Big(\1n H(x,\bar y(x),\bar\psi(x),\1n\nabla\bar y(x),\1n\nabla\bar\psi(x),\1n\bar u(x))\1n-\1n H(x,\bar y(x),\1n\bar\psi(x),\1n\nabla\bar y(x),\1n\nabla\bar\psi(x),\1n u(x))\1n\Big){\ell(x)^\top\2n A(x,\1n\bar u(x))\ell(x)\over\ell(x)^\top\2n A(x,\1n u(x))\ell(x)}\\
\ns\ds\qq+\Big(H_y(x,\bar y(x),\bar\psi(x),\nabla\bar y(x),\nabla\bar\psi(x),\bar u(x))- H_y(x,\bar y(x),\bar\psi(x),\nabla\bar y(x),\nabla\bar\psi(x),u(x))\Big)Y(x)\\
\ns\ds\qq-{1\over 2}H_{yy}(x,\bar y(x),\bar\psi(x),\nabla\bar y(x),\nabla \bar\psi(x),\bar u(x))|Y(x)|^2\\
\ns\ds\qq+\lan\big(A(x,u(x))-A(x,\bar u(x)\big)\nabla\bar\psi(x),\nabla Y(x)\ran\Big\}dx\ges0,\ea\ee
where $\bar\psi(\cd)$ is the weak solution to the adjoint equation \rf{adjoint}, $H(\cd)$ is the Hamiltonian defined by \rf{H}, and $Y(\cd)$ is the weak solution to the following {\it variational equation}:
\bel{variation}\left\{\2n\ba{ll}
\ns\ds-\Div\big(A(x,\bar u(x))\nabla Y(x)\big)=f_y(x,\bar y(x),\bar u(x))\,Y(x)+\Div\big(
\Th(x)\nabla\bar y(x)\big) \\
\ns\ds\qq\qq\qq\qq\qq\qq\q+f(x,\bar y(x),u(x))-f(x,\bar y(x),\bar u(x)),\qq\eqin\O,\\
\ds Y\big|_{\pa\O}=0,\ea\right.\ee
with
\bel{Th}\ba{ll}
\ns\ds\Th(x)=A(x,u(x))-A(x,\bar u(x))-{\big[A(x,u(x))-A(x,\bar u(x))\big]\ell(x)\ell(x)^\top\big[A(x,u(x))-A(x,\bar u(x))\big]\over
\ell(x)^\top A(x,u(x))\ell(x)}\,.\ea\ee

\et

The proof of the above theorem will be carried out in Section 5. The following is a result concerning the partially singular (instead of partially weakly singular) optimal controls.

\bc{T205} \sl Let {\rm(S1)--(S5)} hold, and $(\by(\cd),\bu(\cd))$ be an optimal pair of Problem {\rm(C)}, with $\bu(\cd)$ being partially singular at $u(\cd)$. Then
\bel{second2}\ba{ll}
\ns\ds\int_\O\[\(H_y(x,\by(x),\bar\psi(x),\nabla\by(x),\nabla\bar\psi(x),\bu(x))- H_y(x,\by(x),\bar\psi(x),\nabla\by(x),\nabla\bar\psi(x),u(x))\)Y(x)\\
\ns\ds\qq-{1\over2}H_{yy}(x,\by(x),\bar\psi(x),\nabla\by(x),\nabla\bar\psi(x),\bu(x))|Y(x)|^2\\
\ns\ds\qq+\lan\big(A(x,u(x))-A(x,\bu(x)\big)\nabla\bar\psi(x),\nabla Y(x)\ran\]dx\ges0,\ea\ee
where $\bar\psi(\cd)$ and $H(\cd)$ are the same as those in Theorem \ref{T204}, and $Y(\cd)$ is the weak solution to the following variational equation:
\bel{variation2}\left\{\2n\ba{ll}
\ns\ds-\1n\Div\big(A(x,\bu(x))\nabla Y(x)\big)\1n=\2n f_y(x,\by(x),\bu(x))Y(x)\1n+\2n\Div \big(
A(x,u(x))\nabla\by(x)\big)+f(x,\by(x),u(x)),~\eqin  \O,\\
\ns\ds Y |_{\pa\O}=0. \ea\right.\ee

\ec

\it Proof. \rm By Proposition \ref{P203}, we know that since $\bu(\cd)$ is singular at $u(\cd)$, \rf{(A-A)=0} holds for some $\ell(\cd)\in\cL$. Then
$$\Div\big(\Th(s)\nabla\bar y(x)\big)=\Div\([A(x,u(x))-A(x,\bar u(x))]\nabla\bar y(x)\)
=\Div\(A(x,u(x))\nabla\bar y(x)\)+f(x,\bar y(x),\bar u(x)).$$
Hence, \rf{variation} becomes \rf{variation2}, and \rf{second} becomes \rf{second2} (making use the singularity of $\bu(\cd)$ at $u(\cd)$, see \rf{singular}). Therefore, our conclusion follows. \endpf

\vskip-0.5cm

The following gives the situation that the leading term does not contain the control, whose proof is pretty straightforward.

\bc{T404} \sl Let {\rm(S1)--(S5)} hold with $A(x,v)\equiv A(x)$ independent of $v$. Let $(\by(\cd),\bu(\cd))$ be an optimal pair of Problem {\rm(C)}, for which $\bu(\cd)$ is singular at
$u(\cd)\in\cU_0(\bar u(\cd))\setminus\{\bu(\cd)\}$. Then
\bel{second3}\ba{ll}
\ns\ds\int_\O\[\(\cH_y(x,\by(x),\bar\psi(x),\bu(x))-\cH_y(x,\by(x),\bar\psi(x),u(x))\)Y(x)dx\\
\ns\ds\qq-{1\over 2}\cH_{yy}(x,\by(x),\bar\psi(x),\nabla\by(x),\nabla\bar\psi(x),\bu(x))|Y(x)|^2\]dx\ges0,\ea\ee
where
\bel{H2}\cH(x,y,\psi,v)=\lan\psi,f(x,y,v)\ran-f^0(x,y,v),\q(x,y,\psi,v)\in\dbR^n\times\dbR\times \dbR\times U,\ee
and $\bar\psi(\cd)$ and $Y(\cd)$ are the weak solutions to the following adjoint equation and variational equation, respectively:
\bel{adjoint2}\left\{\2n\ba{ll}
\ns\ds-\Div\big(A(x)\nabla\bar\psi(x)\big)=f_y(x,\by(x),\bu(x))\,\bar\psi(x)
-f^0_y(x,\by(x),\bu(x)),\q\eqin  \O,\\
\ns\ds\bar\psi\big|_{\,\pa\O}=0\ea\right.\ee
\bel{variation3}\left\{\2n\ba{ll}
\ns\ds-\Div\big(A(x)\nabla Y(x)\big)=f_y(x,\by(x),\bu(x))\,Y(x)+f(x,\by(x),u(x))-f(x,\by(x),\bu(x)),\q\hb{in }  \O,\\
\ns\ds Y |_{\pa\O}=0.\ea\right.\ee

\ec

\section{The First-Order Necessary Condition Revisited}

In this section, we briefly recall the proof of Theorem \ref{T201}, from which we will find a correct direction approaching the second-order necessary condition for the optimal control. To this end, we first recall the following lemma (\cite{Allaire 1992}).

\ms

\bl{T301} \sl Let $G:\O\times\dbR^n\to\dbS^n_+$ be measurable. Assume that

\ms

{\rm(i)} $z\mapsto G(x,z)$ is $[0,1]^n$--periodic;

\ms

{\rm(ii)} There exists two constants $\L>\l>0$, such that
\bel{E301}\l|\xi|^2\les\lan G(x,z)\xi,\xi\ran\les\L |\xi|^2,\qq\forall(x,z)\in
\O\times\dbR^n,~\xi\in\dbR^n;\ee

{\rm(iii)} The following holds:
\bel{E302}\lim_{\e\to0^+}\int_\O\Big|G\big(x,{x\over\e}\big)\Big|^2\,dx=\int_\O\int_{[0,1]^n}
|G(x,s)|^2\,ds\,dx.\ee
Let $g\in H^{-1}(\O)$ and $y_\e(\cd)$ be the solution of
\bel{E303}\left\{\2n\ba{ll}
\ns\ds-\Div\(G\big(x,{x\over\e}\big)\nabla y_\e(x)\)=g,\qq\eqin\O,\\
\ns\ds y_\e\big|_{\,\pa\O}=0.\ea\right.\ee
Then $y_\e(\cd)$ converges weakly to $y(\cd)$ in $H^1_0(\O)$ where
$y(\cd)$ solves
\bel{E304}\left\{\2n\ba{ll}
%215
\ns\ds-\Div\big(\h G(x)\nabla y(x)\big)=g,\qq\eqin\O,\\
\ns\ds y\big|_{\,\pa\O}=0,\ea\right.\ee
and $\h G(x)\equiv\(\h G_{ij}(x)\)$ is given by
\bel{E305}\h G_{ij}(x)=\int_{[0,1]^n}\lan G(x,z)\nabla_z[\phi_i(x,z)+z_i\big],\nabla_z[\phi_j(x,z)+z_j]\ran dz,\qq1\les i,j\les n,\ee
%216
with $\phi_i(x,\cd)\in W^{1,2}_\#([0,1]^n;\dbR^n)/\dbR^n$ being the unique solution of
\bel{E306}
%217
-\nabla_z\cd\(G(x,z)\nabla_z[\phi_i(x,z)+z_i]\)=0,\qq1\les i\les n.\ee

\el

Observe that
$$\ba{ll}
\ns\ds\int_{[0,1]^n}\lan G(x,z)\nabla_z[\phi_i(x,z)+z_i],\nabla_z[\phi_j(x,z)+z_j]\ran dz\\
\ns\ds=\int_{[0,1]^n}\(\lan\nabla_z\phi_i(x,z),G(x,z)\big[\nabla_z\phi_j(x,z)+z_j\big]\ran
+\lan G(x,z)e_i,\nabla_z\phi_j(x)+e_j\ran\)dz\\
\ns\ds=\int_{[0,1]^n}\(\lan G(x,z)e_i,e_j\ran+\lan G(x,z)e_i,\nabla_z\phi_j(x,z)\ran\)dz\\
\ns\ds=\int_{[0,1]^n}\(G_{ij}(x,z)+\big(G(x,z)\nabla_z\phi(x,z)^\top\big)_{ij}\)dz.\ea$$
Hence,
\bel{E305*}\h G(x)=\int_{[0,1]^n}G(x,z)\(I+\nabla_z\phi(x,z)^\top\)dz.\ee
Also, \rf{E306} can be written as
\bel{E306*}-\nabla_z^\top\1n\[G(x,z)(I+\nabla_z\phi^\top)\]=0.\ee

Note that in general $G(x,{x\over\e})$ does not necessarily converge strongly in $L^2(\O)$ (as $\e\da0$). Therefore, the above lemma is by no means trivial or obvious. On the other hand, the following result is much easier, which will also be used later, for different situations.

\bl{T302} \sl Let $\e>0$ and $G_\e(\cd)\in L^\infty(\O;\dbS^n_+)$. Assume that there exist two constants $\L>\l>0$, such that
$$\l|\xi|^2\les\lan G_\e(x)\xi,\xi\ran\les\L|\xi|^2, \qq\forall x\in\O,\xi\in\dbR^n;$$
and $G_\e(\cd)$ converges to $G(\cd)$ strongly in $L^2(\O)$. Let $g\in H^{-1}(\O)$ and $y_\e(\cd)$ be the solution of
$$\left\{\2n\ba{ll}
\ns\ds-\Div\(G_\e(x)\nabla y_\e(x)\)=g,\qq\eqin  \O,\\
\ns\ds y_\e\big|_{\,\pa\O}=0.\ea\right.$$
Then $y_\e(\cd)$ converges strongly to $y(\cd)$ in $H^1_0(\O)$, as $\e\da0$, where $y(\cd)$ solves
$$\left\{\2n\ba{ll}
\ns\ds-\Div\(G(x)\nabla y(x)\)=g,\qq\eqin \O,\\
\ns\ds y\big|_{\,\pa\O}=0.\ea\right.$$

\el

In the above lemma, one can prove easily that $y_\e(\cd)$ converges weakly to $y(\cd)$ in $H^1_0(\O)$, as $\e\da0$. The strong convergence follows from
$$\ba{ll}
\ns\ds\lim_{\e\da0}\int_\O\lan G_\e(x)\big(\nabla y_\e(x)-\nabla y(x)\big),\nabla
y_\e(x)-\nabla y(x)\ran\,dx\\
\ns\ds=\lim_{\e\da0}\int_\O\lan G_\e(x)\nabla y_\e(x),\nabla
y_\e(x)-\nabla y(x)\ran\,dx=\lim_{\e\da0}\int_\O(y_\e(x)-y(x))g\,dx=0.\ea$$

%Note that between Lemmas \ref{T301} and \ref{T302}, none of them covers the other.

\ms

We now recall the proof of Theorem \ref{T201} (see \cite{Lou-Yong 2009} for technical details). Let $\bar u(\cd)\in\cU$ be an optimal control and $u(\cd)\in\cU$ be an arbitrary fixed control. Pick any $\m\in S^{n-1}$. Define a two-parameter spike variation $u^{\a,\e}(\cd\,;\m)$ of the control $\bar u(\cd)$ associated with $u(\cd)$ and $\m$ as follows:
\bel{u^(a,e)}
%220
u^{\a,\e}(x;\m)=\left\{\2n\ba{ll}
\ns\ds u(x),\qq\hb{if }
\Big\{{\lan x,\m\ran\over\e}\Big\}\in [0,\a),\\
\ns\ds\bar u(x),\qq\hb{if }\Big\{{\lan x,\m\ran\over\e}\Big\}\in[\a,1),\ea\right.\ee
where $\{a\}\equiv a-[a]$ denotes the decimal part of the real number $a$. Then $u^{\a,\e}(\cd\,;\m)\in\cU$. Here, the dependence on $\m$ is emphasized. We should keep in mind that $u^{\a,\e}(\cd\,;\m)$ also depends on the selected control $u(\cd)$ (which is fixed). Let $y^{\a,\e}(\cd\,;\m)=y(\cd;u^{\a,\e}(\cd\,;\m))$ be the state corresponding to the control $u^{\a,\e}(\cd\,;\m)$. Then by Lemma \ref{T301}, as $\e\da0$, $y^{\a,\e}(\cd\,;\m)$ converges to $y^\a(\cd\,;\m)$, weakly in $H^1_0(\O)$ and strongly in $L^2(\O)$, where $y^\a(\cd\,;\m)$ solves the following PDE, which is called a {\it relaxed state equation}:
\bel{state-y^a(m)}\left\{\2n\ba{ll}
%221
\ns\ds-\Div\big(A^\a(x;\m)\nabla y^\a(x;\m)\big)=(1-\a)f(x,y^\a(x;\m),\bar u(x))+\a f(x,y^\a(x;\m),u(x)),\qq\eqin\O,\\
\ns\ds y^\a(x;\m)=0,\qq\eqin\pa\O,\ea\right.\ee
with
\bel{A^a(m)}\ba{ll}
%222
\ns\ds A^\a(x;\m)=\a A(x,u(x))+(1-\a)A(x,\bar u(x))\\
\ns\ds\qq\qq\q-{\a(1-\a)\big[A(x,u(x))-A(x,\bar u(x))\big]\m \m^\top
\big[A(x,u(x))-A(x,\bar u(x))\big]\over(1-\a)\m^\top A(x,u(x))\m+\a
\m^\top A (x,\bar u(x))\m}.\ea\ee
Define
$$Y^\a(\cd\,;\m)={y^\a(\cd\,;\m)-\bar y(\cd)\over\a}.$$
Then, as $\a\da0$, $Y^\a(\cd\,;\m)$ converges to $Y(\cd\,;\m)$ weakly in $H^1_0(\O)$ and strongly in $L^2(\O)$, where $Y(\cd\,;\m)$ is the weak solution to the following:
\bel{variation-m}\left\{\2n\ba{ll}
\ns\ds-\Div\big(A(x,\bar u(x))\nabla Y(x;\m)\big)=f_y(x,\bar y(x),\bar u(x))\,Y(x;\m)+\Div\big(\Th(x;\m)\nabla\bar y(x)\big) \\
\ns\ds\qq\qq\qq\qq\qq\qq\q+f(x,\bar y(x),u(x))-f(x,\bar y(x),\bar u(x)),\qq\eqin\O,\\
\ds Y(x;\m)=0,\qq\qq\eqin\pa\O,\ea\right.\ee
with
\bel{Th(m)}\ba{ll}
\ns\ds\Th(x;\m)=A(x,u(x))-A(x,\bar u(x))-{\big[A(x,u(x))-A(x,\bar u(x))\big]\m\m^\top\big[A(x,u(x))-A(x,\bar u(x))\big]\over
\m^\top A(x,u(x))\m}\,.\ea\ee
Consequently, as $\a\da0$, $y^\a(\cd\,;\m)$ converges to $\bar y(\cd)$ weakly in $H_0^1(\O)$ and strongly in $L^2(\O)$. On the other hand, by the convergence of $y^{\a,\e}(\cd\,;\m)\to y^\a(\cd\,;\m)$ (strongly in $L^2(\O)$, as $\e\da0$), one has
\bel{lim(e)J(m)}\lim_{\e\da0}J(u^{\a,\e}(\cd\,;\m))=J^\a(u(\cd),\m)\equiv\a\int_\O f^0(x,y^\a(x;\m),u(x))\,dx+(1-\a)\int_\O f^0(x,y^\a(x;\m),\bar u(x))\,dx.\ee
Further, (suppressing $x$)
\bel{J^a(m)-J}\ba{ll}
\ns\ds J^\a(u(\cd),\m)-J(\bu(\cd))=\int_\O\[\a\(f^0(y^\a,u)-f^0(y^\a,\bu)\)
+f^0(y^\a,\bu)-f^0(\by,\bu)\]dx\\
\ns\ds\qq\qq\qq\qq\q=\a\int_\O\[\(f^0(y^\a,u)-f^0(y^\a,\bu)\)+\(\int_0^1f^0_y(\by+\a tY^\a,\bu)
dt\)Y^\a\]dx.\ea\ee
Thus,
\bel{lim(a)(J^a(m)-J)/a}\lim_{\a\da0}{J^\a(u(\cd),\m)-J(\bar u(\cd))\over\a}=\int_\O\[f^0(x,\bar y(x),u(x))-f^0(x,\bar y(x),\bar u(x))+f^0_y(x,\bar y(x),\bar u(x))Y(x;\m)\]dx.\ee
Consequently,
\bel{J^(a,e)-J=}\ba{ll}
\ns\ds J\big(u^{\a,\e}(\cd\,;\m)\big)-J(\bar u(\cd))=J^\a(u(\cd),\m)-J(\bar u(\cd))+r_\e\\
\ns\ds=\a\int_\O\[f^0(x,\bar y(x),u(x))-f^0(x,\bar y(x),\bar u(x))
+f^0_y(x,\bar y(x),\bar u(x))Y(x;\m)\]dx+\a\rho_\a+r_\e,\ea\ee
where
$$\lim_{\e\da0}r_\e=0,\qq\lim_{\a\da0}\rho_\a=0.$$
By the duality, we obtain
\bel{J^(a,e)=J+aJ^1+}\ba{ll}
\ns\ds J(u^{\a,\e}(\cd\,;\m))=J(\bar u(\cd))\\
\ns\ds\qq\q+\a\2n\int_\O\[H(x,\bar y(x),\bar\psi(x),
\nabla\bar y(x),\nabla\bar\psi(x),\bar u(x))
-H(x,\bar y(x),\bar\psi(x),\nabla\bar y(x),\nabla\bar\psi(x),u(x))\\
\ns\ds\qq\q-{\lan\big[A(x,\bar u(x))-A(x,u(x))\big]\nabla\bar y(x),\m\ran\,
\lan\big[A(x,\bar u(x))-A(x,u(x))\big]\nabla\bar\psi(x),\m\ran\over\lan A(x,u(x))\m,\m\ran}\]dx+\a \rho_\a+r_\e\\
\ns\ds\equiv J(\bu(\cd))+\a J^1(\bu(\cd);u(\cd),\m)+\a\rho_\a+r_\e.\ea\ee
By the optimality of $\bar u(\cd)$, the above leads to the following:
$$J^1(\bu(\cd);u(\cd),\m)\ges0,$$
which implies
\bel{H-H=(A-A)(m)>0}\ba{ll}
%226
\ns\ds H\big(x,\bar y(x),\bar\psi(x),\nabla\bar y(x),\nabla\bar\psi(x),\bar u(x)\big)
-H\big(x,\bar y(x),\bar\psi(x),\nabla\bar y(x),\nabla\bar\psi(x),v\big)\\ [2mm]
\ns\ds-{\lan\big[A(x,\bar u(x))-A(x,v)\big]\nabla\bar y(x),\m\ran\,
\lan\big[A(x,\bar u(x))-A(x,v)\big]\nabla\bar\psi(x),\m\ran\over\lan A(x,v)
\m,\m\ran}\ges0,\\
\ns\ds\qq\qq\qq\qq\qq\qq\qq\qq\qq\qq\forall v\in U,\q\m\in S^{n-1},\q\ae~x\in\O.\ea\ee
Hence, \rf{H-H>0} follows, proving the first-order necessary condition.

\ms

Note that in the above result, $\m\in S^{n-1}$ is a given fixed direction. Whereas, when an optimal control $\bu(\cd)$ is weakly singular at $(u(\cd),\ell(\cd))\in\cV_0(\bu(\cd))$, $\ell(\cd)\in\cL$ might not be a fixed $\m$. Therefore, we need to extend \rf{J^(a,e)=J+aJ^1+}, allowing $\m$ to be replaced by $\ell(\cd)\in\cL$. More precisely, we hope to have the following result.

\bp{T303} \sl Let {\rm(S1)--(S4)} hold and $(\by(\cd),\bu(\cd))$ be an optimal pair of Problem {\rm(C)}. Let $(u(\cd),\ell(\cd))\in\cU\times\cL$, and let $y^\a(\cd\,;\ell(\cd))$ be the weak solution to the following equation:
\bel{state-y^a(ell)}\left\{\2n\ba{ll}
\ns\ds-\Div\(A^\a(x;\ell(\cd))\nabla y^\a(x;\ell(\cd))\)=(1-\a)f\big(x,y^\a(x;\ell(\cd)),\bar u(x)\big)\1n+\1n\a f\big(x,y^\a(x;\ell(\cd)),u(x)\big),\q\eqin\O,\\
\ns\ds y^\a(x;\ell(\cd))=0,\qq\qq\eqin\pa\O.\ea\right.\ee
with
\bel{A^a(ell)}\ba{ll}
\ns\ds A^\a(x;\ell(\cd))=\a A(x,u(x))+(1-\a)A(x,\bu(x))\\
\ns\ds\qq\qq\qq-{\a(1-\a)\big[A(x,u(x))-A(x,\bu(x))\big]\ell(x) \ell(x)^\top
\big[A(x,u(x))-A(x,\bu(x))\big] \over (1-\a) \ell(x)^\top
A(x,u(x))\ell(x)+\a \ell(x)^\top A (x,\bu(x))\ell(x)}.\ea\ee
Define
\bel{J^a(u,ell)}J^\a\big(u(\cd);\ell(\cd)\big)=\a\int_\O f^0\big(s,y^\a(x;\ell(\cd)),u(x)
\big)dx+(1-\a)\int_\O f^0\big(x,y^\a(x;\ell(\cd)),\bar u(x)\big)dx.\ee
Then
\bel{J^a(u,ell)=J+aJ^1+}\ba{ll}
\ns\ds J^\a(u(\cd),\ell(\cd))=J(\bar u(\cd))\\
\ns\ds\q+\a\int_\O\[H(x,\bar y(x),\bar\psi(x),
\nabla\bar y(x),\nabla\bar\psi(x),\bar u(x))
-H(x,\bar y(x),\bar\psi(x),\nabla\bar y(x),\nabla\bar\psi(x),u(x))\\
\ns\ds\q-{\lan\big[A(x,\bar u(x))-A(x,u(x))\big]\nabla\bar y(x),\ell(x)\ran\,
\lan\big[A(x,\bar u(x))-A(x,u(x))\big]\nabla\bar\psi(x),\ell(x)\ran\over\lan A(x,u(x))\ell(x),\ell(x)\ran}\]dx\1n+\1n\a\rho_\a,\ea\ee
and
\bel{J^a(u,ell)-J>0}J^\a(u(\cd),\ell(\cd))-J(\bar u(\cd))\ges0.\ee

\ep

If the above result holds true, then in the case that $\bu(\cd)$ is weakly singular at $(u(\cd),\ell(\cd))$, the above \rf{J^a(u,ell)=J+aJ^1+} will become the following in which the first order term disappears
$$J^\a(u(\cd),\ell(\cd))=J(\bu(\cd))+\a\rho_\a.$$
To further characterize the optimal control, the second-order necessary condition will be needed.

\ms

To prove Proposition \ref{T303}, it is natural to try a modification of \rf{u^(a,e)} as follows:
\bel{u^(a,e)(ell)}u^{\a,\e}(x)=\left\{\2n\ba{ll}
\ns\ds u(x),\qq\hb{if }\Big\{{\lan x,\ell(x)\ran\over\e}\Big\}\in [0,\a),\\
\ns\ds\bu(x),\qq\hb{if }\Big\{{\lan x,\ell(x)\ran\over\e}\Big\}\in[\a,1),\ea\right.\ee
and wish that $y^{\a,\e}(\cd)=y(\cd;u^{\a,\e}(\cd))$ converges to $y^\a(\cd)$ weakly in
$H^1_0(\O)$ (as $\e\da0$) with $y^\a(\cd)$ being the weak solution to \rf{state-y^a(ell)}. However, \eqref{u^(a,e)(ell)} does not work as we expected. For example, for $n=2$, let
$$\ell(x_1,x_2)={(-x_2,x_1)\over\sqrt{x_1^2+x_2^2}},\qq\forall x\ne0.$$
Then \eqref{u^(a,e)(ell)} implies
$$u^{\a,\ve}(x)=u(x),\qq\forall x\in\O.$$
Thus, in this case, $u^{\a,\e}(\cd)$ is not a proper perturbation of $\bu(\cd)$ that we
expected since for any proper metric $\wt\rho$ on $\cU$,
$$\wt\rho(u^{\e,\a}(\cd),\bu(\cd))=\wt\rho(u(\cd),\bu(\cd)),$$
which will not go to zero as $\a,\e\da0$. Nevertheless, in the next section, we will prove that  Proposition \ref{T303} is true, by a different method.

\section{Proof of Proposition \ref{T303}.}

In this section, we will present a proof of Proposition \ref{T303}. Let us begin with the following lemma.

\bl{T401} \sl Let $B_1,B_2\in\dbS^n_+$ and
$$\l|\xi|^2\les\lan B_i\xi,\xi\ran\les\L|\xi|^2,\qq\forall i=1,2,\q\xi\in\dbR^n$$
for some constants $\L\ges\l>0$. Let $\a\in (0,1)$, $\m\in\dbR^n\setminus\{0\}$ and
$$G=\a B_1+(1-\a)B_2-\a(1-\a){\big(B_2-B_1\big)\m\m^\top
\big(B_2-B_1\big)\over\m^\top\big[\a B_2+(1-\a)B_1\big]\m}.$$
Then
\if{
\bel{E416}\l|\xi|^2\les\lan G\xi,\xi\ran\les\L|\xi|^2,\qq\forall\xi\in\dbR^n.\ee
\el
}\fi

\bel{E416} \big(\a B_1^{-1}+(1-\a)B_2^{-1}\big)^{-1}\leq  G\leq \a B_1+(1-\a)B_2.\ee
\el

\it Proof. \rm The lemma is a consequence of Theorem 1.3.14 and Lemma 1.3.32 in \cite{Allaire 2002}. Here we give a direct proof of it.
Let
$$C\equiv\a B_2+(1-\a)B_1=B_1+\a(B_2-B_1)=B_2-(1-\a)(B_2-B_1).$$
Then $C\in\dbS_+^n$, and
$$C-B_1=\a(B_2-B_1),\qq C-B_2=-(1-\a)(B_2-B_1).$$
Thus,
\bel{E417}\ba{ll}
\ns\ds\a B_1+(1-\a)B_2-\a(1-\a)\big(B_2-B_1\big)C^{-1}\big(B_2-B_1\big)\\
\ns\ds=\a B_1+(1-\a)B_2+(C-B_1)C^{-1}(C-B_2)\\
\ns\ds=\a B_1+(1-\a)B_2+C-B_1-B_2+B_1C^{-1} B_2\\
\ns\ds=B_1[\a B_2+(1-\a)B_1]^{-1} B_2=[\a B^{-1}_1+(1-\a)B_2^{-1}]^{-1}.\ea\ee
Consequently, we get \eqref{E416} since
\bel{E418}\ba{ll}
\ns\ds 0\leq {\m\m^\top
\over\m^\top\big[\a B_2+(1-\a)B_1\big]\m}={ \m  \m ^\top  \over\m^\top C\m}=C^{-{1\over  2}}{C^{1\over  2}\m\over |C^{1\over  2}\m |} \Big( {C^{1\over  2}\m\over |C^{1\over  2}\m |}\Big)^\top C^{-{1\over  2}}
\leq C^{-1},\ea\ee
proving our conclusion. \endpf

\vskip-0.5cm

The following lemma will play an interesting role blow.

\bl{T402} \sl Let $\n=(\n_1,\cds,\n_n)^\top\in\dbZ^n\setminus\{0\}$ (where $\dbZ$ is the set of all integers). Then for any $\a\in(0,1)$,
$$\int_{[0,1]^n}\chi_{[0,\a)}(\{\lan\n,z\ran\})dz=\a.$$

\el

\it Proof. \rm Recall that $\{a\}=a-[a]$ is the decimal part of the real number $a$. Note that
$$\lan\n,z\ran=\sum_{k=1}^n\n_kz_k,$$
If some of integers $\n_k$ are zero, we could drop the corresponding terms and reduce the dimension of $z$. Thus, we assume all $\n_k$ are non-zero. Also, if some $\n_k<0$, we may replace corresponding $z_k$ by $(1-z_k)$. Therefore, we may let all $\n_k>0$. Next, we observe the following (noting the $[0,1]^n$-periodicity of the maps $\ds z\mapsto\big\{\lan\n,z\ran\big\}$):
$$\ba{ll}
\ns\ds\int_{[0,1]^n}\chi_{[0,\a)}\big(\{\lan\n,z\ran\}\big)dz=\int_0^1dz_1\int_0^1dz_2\cds\int_0^1
\chi_{[0,\a)}\(\Big\{\sum_{k=1}^n\n_kz_k\Big\}\)dz_n\\
\ns\ds=\n_1\n_2\cds\n_n\int_0^{1\over\n_1}dz_1\int_0^{1\over\n_2}ds_2\cds\int_0^{1\over\n_n}\chi_{[0,\a)}
\(\Big\{\sum_{k=1}^n\n_kz_k\Big\}\)dz_n\\
\ns\ds=\int_0^1dz_1\int_0^1dz_2\cds\int_0^1\chi_{[0,\a)}\(\Big\{\sum_{k=1}^nz_k\Big\}\)dz_n.\ea$$
Hence, we need to prove the following:
$$\int_0^1dz_1\int_0^1dz_2\cds\int_0^1\chi_{[0,\a)}\(\Big\{\sum_{k=1}^nz_k\Big\}\)dz_n=\a.$$
Let us use induction. For $n=1$, the above is clearly true. Suppose the above holds for $n-1$. Then, for the $n$-dimensional case, we observe the following: For $z_1\in[0,1)$,
$$\Big\{\sum_{k=1}^nz_k\Big\}=\Big\{z_1+\sum_{k=2}^nz_k\Big\}=\left\{\2n\ba{ll}
\ns\ds z_1+\Big\{\sum_{k=2}^nz_k\Big\},\qq\qq\hb{if }~0\les z_1+\Big\{\sum_{k=2}^nz_k\Big\}<1,\\
\ns\ds z_1+\Big\{\sum_{k=2}^nz_k\Big\}-1,\qq\hb{if }~1\les z_1+\Big\{\sum_{k=2}^nz_k\Big\}<2.\ea\right.$$
Then for $z_1\in[0,\a)$, the following holds
$$\chi_{[0,\a)}\Big(\Big\{z_1+\sum_{k=2}^nz_k\Big\}\)=1,$$
if and only if either
$$0\les z_1+\Big\{\sum_{k=2}^nz_k\Big\}<\a\qq\hb{and}\qq0\les z_1+\Big\{\sum_{k=2}^nz_k\Big\}<1.$$
or
$$0\les z_1+\Big\{\sum_{k=2}^nz_k\Big\}-1<\a\qq\hb{and}\qq1\les z_1+\Big\{\sum_{k=2}^nz_k\Big\}<2.$$
That is, either
$$0\les\Big\{\sum_{k=2}^nz_k\Big\}<\a-z_1,$$
or
$$1-z_1\les\Big\{\sum_{k=2}^nz_k\Big\}\les(\a+1-z_1)\land1=1.$$
Note that the above two cases are mutually exclusive (since $\a-z_1<1-z_1$). On the other hand, for $z_1\in[\a,1)$, if $\ds0\les z_1+\Big\{\sum_{k=2}^nz_k\Big\}<1$,
then
$$\Big\{z_1+\sum_{k=2}^nz_k\Big\}=z_1+\Big\{\sum_{k=2}^nz_k\Big\}\ges\a\q\Ra\q
\chi_{[0,\a)}\(z_1+\Big\{\sum_{k=2}^nz_k\Big\}\)=0.$$
Thus the following holds:
$$\chi_{[0,\a)}\Big(\Big\{z_1+\sum_{k=2}^nz_k\Big\}\)=1,$$
if and only if
$$0\les z_1+\Big\{\sum_{k=2}^nz_k\Big\}-1<\a,\q1\les z_1+\Big\{\sum_{k=2}^nz_k\Big\}<2.$$
That is,
$$0\les1-z_1\les\Big\{\sum_{k=2}^nz_k\Big\}\les\a+1-z_1\les1.$$
Hence, by induction hypothesis,
$$\ba{ll}
\ns\ds\int_{[0,1]^n}\chi_{[0,\a)}\(\Big\{z_1+\sum_{k=2}^nz_k\Big\}\)dz
=\int_0^\a dz_1\int_{[0,1]^{n-1}}\chi_{[0,\a-z_1)}\(\Big\{\sum_{k=2}^nz_k\Big\}\)
dz_2\cds dz_n\\
\ns\ds\qq\qq\qq\qq\qq\qq\qq+\int_0^\a dz_1\int_{[0,1]^{n-1}}\chi_{[1-z_1,1)}\(\Big\{\sum_{k=2}^n
z_k\Big\}\)dz_2\cds dz_n\\
\ns\ds\qq\qq\qq\qq\qq\qq\qq+\int_\a^1dz_2\int_{[0,1]^{n-1}}\chi_{[1-z_1,\a+1-z_1)}\(\Big\{\sum_{k=1}^n
z_k\Big\}\)dz_2\cds dz_n\\
\ns\ds=\int_0^\a(\a-z_1)dz_1+\int_0^\a z_1dz_1+\int_\a^1\a dz_1=\a.\ea$$
This completes the proof. \endpf

\vskip-0.5cm

The following gives a crucial convergence of the weak solution to the state equation under a suitable perturbation of the leading coefficient.

\bl{T403} \sl Let {\rm(S1)} hold. Let $B(\cd)=\big(b_{ij}(\cd)\big),C(\cd)=\big(c_{ij}(\cd)\big)\in
L^\infty(\O;\dbS^n_+)$ such that
\bel{E401}\ba{ll}
\ns\ds\l|\xi|^2\les\lan B(x)\xi,\xi\ran\les\L|\xi|^2,\\
\ns\ds\l|\xi|^2\les\lan C(x)\xi,\xi\ran\les\L|\xi|^2,\ea\qq\forall
(\xi,x)\in\dbR^n\times\O,\ee
for some $\L\ges\l>0$. Let $\a\in(0,1)$, $\ds\m_k=\({r_{k1}\over p_{k1}},{r_{k2}\over p_{k2}},\ldots,{r_{kn}\over p_{kn}}\)\in\dbQ^n$,
$\dbQ$ is the set of all rational numbers, with all $r_{kj}$ being integers, and $p_{kj}$ being positive integers, $1\les k\les m$, $1\les j\les n$. Let $\ds E_1,E_2,\ldots, E_m$ be mutually disjoint measurable sets such that $\ds\bigcup^m_{k=1}E_k=\O$. Let $h(\cd)\in L^2(\O)$, and
\bel{E402}G(x,z)=\left\{\2n\ba{ll}
\ns\ds B(x),\qq\hb{if }\big\{\lan z,\m_k\ran\big\}\in[0,\a),\\
\ns\ds C(x),\qq\hb{if }\big\{\lan z,\m_k\ran\big\}\in[\a,1),\ea\right.\q x\in E_k,~1\les k\les m.\ee
For $\e>0$, let $y^\e(\cd)\in H^1_0(\O)$ be the weak solution of the following:
\bel{E403}\left\{\ba{ll}
\ns\ds-\Div\(G\big(x,{x\over\e}\big)\nabla
y^\ve(x)\)=h(x),\qq\eqin \O,\\
\ns\ds y^\e\big|_{\,\pa\O}=0.\ea\right.\ee
Then, as $\e\da0$,
\bel{E404}y^\e(\cd)\to y(\cd),\qq\hb{weakly in}\q H^1_0(\O),\ee
where $y(\cd)$ is the weak solution of
\bel{E405}\left\{\2n\ba{ll}
\ns\ds-\Div\big(\h G(x)\nabla y(x)\big)=h(x),\qq\hb{in }\O,\\
\ns\ds y \big|_{\,\pa\O}=0 \ea\right.\ee
with $\h G(\cd)\in L^\infty(\O;\dbS^n_+)$ given by
\bel{E406}\ba{ll}
\ns\ds\h G(x)=\a B(x)+(1-\a)C(x)-{\a(1-\a)\big[B(x)-C(x)\big]\ell(x)\ell(x)^\top \big[B(x)-C(x)\big]
\over(1-\a)\ell(x)^\top B(x)\ell(x)+\a\ell(x)^\top C(x)\ell(x)},\q x\in\O,\ea\ee
and
$$\ell(x)=\sum_{k=1}^m\m_k\chi_{E_k}(x),\qq x\in\O.$$

\el

\it Proof. \rm The proof is essentially inspired by that of Lemma 1.3.32 of \cite{Allaire 2002}.
%We denote
%
%$$G(x,z)=\big(g_{ij}(x,z)\big),\q B(x)=\big(b_{ij}(x)\big),\q C(x)=\big(c_{ij}(x)\big).$$
%
Let $P$ be a common multiple of $p_{kj}$, $k=1,2,\cds,m$,
$j=1,2,\cds,n$. Then one can verify that $G(x,Pz)$ is
$[0,1]^n$--periodic in $z$, for any $\e>0$, $\ds G\(x,{Pz\over\e}\)$
is measurable, and
$$\l|\xi|^2\les\lan G(x,Pz)\xi,\xi\ran\les \L |\xi|^2, \qq\forall
(x,z)\in\O\times\dbR^n,~\xi\in\dbR^n.$$
Moreover, using  Riemann-Lebesgue's Theorem (see Ch. II, Theorem 4.15 in \cite{Zygmund 2002}),
$$\ba{ll}
\ns\ds\lim_{\e\da0}\int_\O\Big|G\(x,{Px\over\e}\)\Big|^2\,dx=
\lim_{\e\da0}\sum^m_{k=1}\int_{E_k}\Big|G\(x,{Px\over\e}\)\Big|^2\,dx\\
\ns\ds=\lim_{\e\da0}\sum^m_{k=1}\int_{E_k}\Big|{B(x)\chi_{[0,\a)}\(\big\{\lan{Px\over
\e},\m_k\ran\big\}\Big)+C(x)\chi_{[\a,1)}\Big(\big\{\lan{Px\over\e},\m_k\ran}\big\}\Big)\Big|^2\,dx\\
\ns\ds=\lim_{\e\da0}\sum^m_{k=1}\int_{E_k}\Big[
|B(x)|^2\,\chi_{[0,\a)}\Big(\big\{\lan{Px\over
\e},\m_k\ran\big\}\Big)+|C(x)|^2\,\chi_{[\a,1)}\Big(\big\{\lan{Px\over
\e},\m_k\ran\big\}\Big)\,\Big] dx\\
\ns\ds=\sum^m_{k=1}\Big[\int_{E_k}|B(x)|^2\, dx\,
\int_{[0,1]^n}\chi_{[0,\a)}\Big(\{\lan{Pz},\m_k\ran\}\Big)\,dz+\int_{E_k}|C(x)|^2\,dx\,
\int_{[0,1]^n}\chi_{[\a,1)}\Big(\{\lan{Pz},\m_k\ran\}\Big)\,dz\Big] \\
%&=&\ds \sum^m_{i=1}\Big[ \a\int_{E_i} \ia{B(x)}^2\, dx
%+(1-\a)\int_{E_i}\ia{C(x)}^2\, dx\Big] \\
%
\ns\ds=\sum^m_{k=1}\int_{E_k}\2n\[\int_{[0,1]^n}\1n\Big|B(x)\,{\huge
\chi}_{[0,\a)}\Big(\{\lan{Pz},\m_k\ran\Big)\1n+\1n
C(x)\,\chi_{[\a,1)}\Big(\{\lan{Pz},\m_k\ran\}\Big)\Big|^2\,dz\]dx
=\2n\int_\O\1n\int_{[0,1]^n}\2n|G(x,Pz)|^2\,dz\,dx.\ea$$
Thus, $(x,z)\mapsto G(x,Pz)$ satisfies conditions of Lemma \ref{T301}. Using Lemma \ref{T301}, we
get that \eqref{E404}--\rf{E405} with
\bel{E407}\h G(x)=\int_{[0,1]^n}G(x,Pz)\[I+\nabla_z\phi(x,z)^\top\]dz,\ee
with $\phi\in W^{1,2}_\#([0,1]^n;\dbR^n)/\dbR^n$ being the unique solution of
\bel{E408}-\nabla_z^\top\1n\[G(x,Pz)\(I+\nabla_z\phi(x,z)^\top\)\]=0.\ee
That is, for $x\in E_k$,
\bel{E409}\nabla_z^\top\Big[\Big(\chi_{[0,\a)}\big(\{\lan Pz,\m_k\ran\}\big)B(x)+
\chi_{[\a,1)}\big(\{\lan Pz,\m_k\ran\}\big)C(x)\Big)\Big(I+\nabla_z\phi(x,z)^\top\)\]=0.\ee
To solve \rf{E409}, we fix $k$, denote $\timu={\m_k\over
|\m_k|}$ and define
$$\ds\f(x,z)=P|\m_k|\phi\big(x,{z\over P|\m_k|}\big),\qq(x,z)\in\dbR^n\times\dbR^n.$$
Then
$$\nabla_z\f(x,z)^\top=\nabla_z\phi\big(x,{z\over P|\m_k|}\big)^\top,$$
which leads to
$$\nabla_z\phi(s,z)^\top=\nabla_z\f\big(x,P|\m_k|z\big)^\top.$$
Using $\f(\cd)$, equation \rf{E409} can be written as
\bel{E409*}\nabla_z^\top\Big[\Big(\chi_{[0,\a)}\big(\{\lan P|\m_k|z,\timu\ran\}\big)B(x)+
\chi_{[\a,1)}\big(\{\lan P|\m_k|z,\timu\ran\}\big)C(x)\Big)\Big(I+\nabla_z\f(x,P|\m_k|z)^\top\Big)\Big]=0,\ee
which is equivalent to the following:
\bel{E410}\nabla_z^\top\[\(\chi_{[0,\a)}\big(\{\lan z,\timu\ran\}\big)B(x)+
\chi_{[\a,1)}\big(\{\lan z,\timu\ran\}\big)C(x)\)\(I+\nabla_z\f(x,z)^\top\)\]=0.\ee
Since $\timu$ is a unit vector, there exists an orthogonal matrix $Q$ such that
$\timu=Q^\top e_1$ with $e_1=(1,0,\ldots,0)^\top\in\dbR^n$. Let $\ds\tiz=(\tiz_1,\tiz_2,\ldots,\tiz_n)^\top=Qz$, and $\wt\f(x,\tiz)=\f(x,Q^\top\wt z)\equiv
\f(x,z)$. Then $\tiz_1=\lan\tiz,e_1\ran=\lan Qz,Q\timu\ran=\lan z,\timu\ran$, and
%
%$$\nabla_z=Q^\top\nabla_{\tiz},$$
%
%Therefore,
%
$$\nabla_z\f(x,z)^\top=Q^\top\nabla_{\tiz}\f(x,Q^\top\tiz)^\top=Q^\top\nabla_\tiz
\wt\f(x,\tiz)^\top.$$
Thus,
\bel{E411}\nabla_{\tiz}^\top\[Q\(\chi_{[0,\a)}\big(\{\tiz_1\}\big)B(x)+
\chi_{[\a,1)}\big(\{\tiz_1\}\big)C(x)\)\(I+Q^\top\nabla_{\tiz}\wt\f(x,\tiz)^\top\)\]=0.\ee
Since $Q\big(\chi_{[0,\a)}\big(\{\tiz_1\}\big)B(x)+\chi_{[\a,1)}\big(\{\tiz_1\}\big)C(x)\big)$, the coefficient of the above equation, is independent of $\tiz_2,\tiz_3,\ldots,\tiz_n$, by the uniqueness, the solution $\wt\f(x,\cd)$ of equation \rf{E411} must be independent of $\tiz_2,\tiz_3,\ldots,\tiz_n$. Thus \rf{E411} further implies
\bel{E412} {\pa\over\pa\tiz_1}\[e_1^\top Q\(\chi_{[0,\a)}\big(\{\tiz_1\}\big)B(x)
+\chi_{[\a,1)}\big(\{\tiz_1\}\big)C(x)\)\(I+\timu
\big({\pa\wt\f\over\pa\tiz_1}(x,\tiz_1e_1)\big)^\top\)\]=0.\ee
Hence, there exists a constant vector $X\in\dbR^n$ such that (note $Q^\top e_1=\timu$)
$$\ba{ll}
\ns\ds X^\top=e_1^\top Q\(\chi_{[0,\a)}\big(\{\tiz_1\}\big)B(x)+\chi_{[\a,1)}\big(\{\tiz_1\}\big)C(x)\)
\[I+\timu\({\pa\wt\f\over\pa\tiz_1}(x,\tiz_1e_1)\)^\top\]\\
\ns\ds\qq=\timu^\top\(\chi_{[0,\a)}\big(\{\tiz_1\}\big)B(x)+\chi_{[\a,1)}\big(\{\tiz_1\}\big)C(x)\)\\
\ns\ds\qq\q+\timu^\top\(\chi_{[0,\a)}\big(\{\tiz_1\}\big)B(x)+\chi_{[\a,1)}\big(\{\tiz_1\}\big)C(x)\)
\timu\({\pa\wt\f\over\pa\tiz_1}\big(x,\tiz_1e_1\big)\)^\top.\ea$$
Consequently,
$$\ba{ll}
\ns\ds{\pa\wt\f\over\pa\tiz_1}(x,\tiz_1e_1)={X-\(\chi_{[0,\a)}(\{\tiz_1\})B(x)+\chi_{[\a,1)}(\{\tiz_1\})C(x)\)
\timu\over\timu^\top\(\chi_{[0,\a)}(\{\tiz_1\})B(x)+\chi_{[\a,1)}(\{\tiz_1\})C(x)\)\timu}\\
\ns\ds\qq\qq\qq={X-B(x)\timu\over\timu^\top B(x)\timu}\chi{[0,\a)}(\{\tiz_1\})+{X-C(x)\timu\over\timu^\top C(x)\timu}\chi_{[\a,1)}(\{\tiz_1\}).\ea$$
Since $\phi(x,\cd)\in W^{1,2}_\#((0,1)^n;\dbR^n)/\dbR^n$, for $\tiz=(P^2|\m_k|^2,\tiz_2,\cds,\tiz_n)^\top$, we have
$$\ba{ll}
\ns\ds\wt\f(x,\tiz)=\wt\f(x,P^2|\m_k|^2e_1)=\f(x,P^2|\m_k|^2Q^\top e_1)=\f(x,P^2|\m_k|^2\timu)\\
\ns\ds\qq\q=P|\m_k|\phi(x,P\mu_k)=P|\m_k|\phi(x,0)=\f(x,0)=\wt\f(x,0).\ea$$
Hence,
$$\ba{ll}
\ns\ds0={\wt\f(x,P^2|\m_k|^2e_1)-\wt\f(x,0)\over P^2|\m_k|^2}\\
\ns\ds\q={1\over P^2|\m_k|^2}\int_0^{P^2|\m_k|^2}\[{X-B(x)\timu\over\timu^\top B(x)\timu}\chi{[0,\a)}(\{\tiz_1\})+{X-C(x)\timu\over\timu^\top C(x)\timu}\chi_{[\a,1)}(\{\tiz_1\})\]d\tiz_1\\
\ns\ds\q=\({\a\over\timu^\top B(x)\timu}+{1-\a\over\timu^\top C(x)\timu}\)X-\({\a B(x)\over\timu^\top B(x)\timu}+{(1-\a)C(x)\over\timu^\top C(x)\timu}\)\timu.\ea$$
This yields
$$X={\a[\timu^\top C(x)\timu]B(x)\timu+(1-\a)[\timu^\top B(x)\timu]C(x)\timu\over\timu^\top\big[\a C(x)+(1-\a)B(x)\big]\timu},$$
and for any $\tiz\in\dbR^n$ with $\tiz_1\in[0,\a)$,
$$\ba{ll}
\ns\ds\wt\f(x,\tiz)=\wt\f(x,\tiz_1e_1)=\wt\f(x,0)+\int_0^{\tiz_1}{\pa\wt\f\over\pa\tiz_1}
(x,\t e_1)d\t\\
\ns\ds\qq\q=\wt\f(x,0)+\int_0^{\tiz_1}{X-B(x)\timu\over\timu^\top B(x)\timu}d\t=\wt\f(x,0)+\tiz_1{X-B(x)\timu\over\timu^\top B(x)\timu}.\ea$$
For $\tiz_1\in[\a,1)$, we have
$$\ba{ll}
\ns\ds\wt\f(x,\tiz)=\wt\f(x,\a e_1)=\wt\f(x,0)+\a{X-B(x)\timu\over\timu^\top B(x)\timu}+\int_\a^{\tiz_1}{\pa\wt\f\over\pa\tiz_1}(x,\t e_1)d\t\\
\ns\ds\qq\q=\wt\f(x,0)+\a{X-B(x)\timu\over\timu^\top B(x)\timu}+\int_0^{\tiz_1}{X-C(x)\timu\over\timu^\top C(x)\timu}
d\t\\
\ns\ds\qq\q=\wt\f(x,0)+\a{X-B(x)\timu\over\timu^\top B(x)\timu}+(\tiz_1-\a){X-C(x)\timu\over\timu^\top C(x)\timu}.\ea$$
As a result,
$$\ba{ll}
\ns\ds\wt\f(x,e_1)-\wt\f(x,0)=\a{X-B(x)\timu\over\timu^\top B(x)\timu}+(1-\a){X-C(x)\timu\over\timu^\top C(x)\timu}\\
\ns\ds=\({\a\over\timu^\top B(x)\timu}+{(1-\a)\over\timu^\top C(x)\timu}\)X-\[{\a B(x)\timu\over\timu^\top B(x)\timu}+{(1-\a)C(x)\timu\over\timu^\top C(x)\timu}\]\\
\ns\ds={\a\timu^\top C(x)\timu+(1-\a)\timu^\top B(x)\timu\over[\timu^\top B(x)\timu][\timu^\top C(x)\timu]}{\a[\timu^\top C(x)\timu]B(x)\timu+(1-\a)[\timu^\top B(x)\timu]C(x)\timu\over\timu^\top\big[\a C(x)+(1-\a)B(x)\big]\timu}\\
\ns\ds\qq-\[{\a B(x)\timu\over\timu^\top B(x)\timu}+{(1-\a)C(x)\timu\over\timu^\top C(x)\timu}\]=0.\ea$$
Thus, $\tiz_1\mapsto\wt\f(x,\tiz_1e_1)$ is 1-periodic.
%Consequently, for any $\tiz\in\dbR^n$,
%
%$$\ba{ll}
%
%\ns\ds\wt\f(x,\tiz)=\wt\f(x,\tiz_1e_1)=\wt\f(x,\{\tiz_1\}e_1)\\
%
%\ns\ds=\wt\f(x,0)+\{\tiz_1\}{X-B(x)\timu\over\timu^\top B(x)\timu}\chi_{[0,\a)}(\{\tiz_1\})
%+\[\a{X-B(x)\timu\over\timu^\top B(x)\timu}+(\{\tiz_1\}-\a){X-C(x)\timu\over\timu^\top C(x)\timu}\]\chi_{[\a,1)}(\{\tiz_1\}).\ea$$
%
On the other hand, for any $z\in[0,1]^n$,
$$\ba{ll}
\ns\ds\phi(x,z)={1\over P|\m_k|}\f\big(x,P|\m_k|z\big)={1\over P|\m_k|}\wt\f\big(x,P|\m_k|Qz\big)\\
\ns\ds\qq\q={1\over P|\m_k|}\wt\f\big(x,P|\m_k|(e_1^\top Qz)e_1\big)={1\over P|\m_k|}
\wt\f\big(x,P|\m_k|(\timu^\top z)e_1\big)={1\over P|\m_k|}
\wt\f\big(x,P|\m_k|(e_1\timu^\top)z\big).\ea$$
Hence,
$$\ba{ll}
\ns\ds\nabla_z\phi(x,z)^\top=\timu e_1^\top\nabla_{\tiz}\wt\f\big(x,(\timu^\top z)
P|\m_k|e_1\big)^\top=\timu{\pa\wt\f\over\pa\tiz_1}\big(x,P(\mu_k^\top z)e_1\big)^\top
=\timu{\pa\wt\f\over\pa\tiz_1}\big(x,\{\lan P\m_k,z\ran\}e_1\big)^\top\\
\ns\ds=\timu{X^\top-\timu^\top\(\chi_{[0,\a)}(\{\lan P\m_k,z\ran\}B(x)
+\chi_{[\a,1)}(\{\lan P\m_k,z\ran\})C(x)\)
\over\timu^\top\(\chi_{[0,\a)}(\{\lan P\m_k,z\ran\})B(x)
+\chi_{[\a,1)}(\{\lan P\m_k,z\ran\})C(x)\)\timu}\\
\ns\ds={\timu[X^\top-\timu^\top B(x)]\over\timu^\top B(x)\timu}\chi_{[0,\a)}(\{\lan P\m_k,z\ran\})+{\timu[X^\top-\timu^\top C(x)]\over\timu^\top C(x)\timu}\chi_{[\a,1)}(\{\lan P\m_k,z\ran\}).\ea$$
Therefore, making use of Lemma \ref{T402}, one has
$$\ba{ll}
\ns\ds\h G(x)=\int_{[0,1]^n}G(x,Pz)\[I+\nabla_z\phi(x,z)^\top\]dz\\
\ns\ds\qq=\int_{[0,1]^n}\[B(x)\chi_{[0,\a)}(\{\lan P\m_k,z\ran\})+C(x)\chi_{[\a,1)}(\{\lan P\m_k,z\ran\})\]\\
\ns\ds\qq\qq\qq\cd\[I+{\timu[X^\top-\timu^\top\top B(x)]\over\timu^\top B(x)\timu}
\chi_{[0,\a)}(\{\lan P\m_k,z\ran\})+{\timu[X^\top-\timu^\top C(x)]\over\timu^\top C(x)
\timu}\chi_{[\a,1)}(\{\lan P\m_k,z\ran\})\]dz\\
\ns\ds\qq=\int_{[0,1]^n}\[\(B(x)+{B(x)\timu[X^\top-\timu^\top B(x)]\over\timu^\top B(x)\timu}\)\chi_{[0,\a)}(\lan\{P\m_k,z\ran\})\\
\ns\ds\qq\qq\qq\q+\(C(x)+{C(x)\timu[X^\top-\timu^\top C(x)]\over\timu^\top C(x)\timu}\)\chi_{[\a,1)}(\lan\{P\m_k,z\ran\})\]dz\\
\ns\ds\qq=\a B(x)+(1-\a)C(x)+\a{B(x)\timu[X^\top-\timu^\top B(x)]\over\timu^\top B(x)\timu}
+(1-\a){C(x)\timu[X^\top-\timu^\top C(x)]\over\timu^\top C(x)\timu}.\ea$$
Now, we simplify the expression of $\h G(x)$,  suppressing $x$,
$$\ba{ll}
\ns\ds B\timu[X^\top-\timu^\top B]=B\timu{\a[\timu^\top C\timu]\timu^\top B
+(1-\a)[\timu^\top B\timu]\timu^\top C\over\timu^\top\big[\a C+(1-\a)B\big]\timu}
-B\timu\timu^\top B\\
\ns\ds={\a[\timu^\top C\timu]B\timu\timu^\top B+(1-\a)[\timu^\top B\timu]B\timu\timu^\top C
-\a[\timu^\top C\timu]B\timu\timu^\top B-(1-\a)[\timu^\top B\timu]B\timu\timu^\top B
\over\timu^\top[\a C+(1-\a)B]\timu}\\
\ns\ds={(1-\a)[\timu^\top B\timu]B\timu\timu^\top(C-B)\over\timu^\top[\a C+(1-\a)B]\timu}.\ea$$
Likewise,
$$\ba{ll}
\ns\ds C\timu[X^\top-\timu^\top C]=C\timu{\a[\timu^\top C\timu]\timu^\top B+(1-\a)
[\timu^\top B\timu]\timu^\top C\over\timu^\top\big[\a C+(1-\a)B\big]\timu}-C\timu\timu^\top C\\
\ns\ds={\a[\timu^\top C\timu]C\timu\timu^\top B+(1-\a)[\timu^\top B\timu]C\timu\timu^\top C
-\a[\timu^\top C\timu]C\timu\timu^\top C-(1-\a)[\timu^\top B\timu]C\timu\timu^\top C
\over\timu^\top[\a C+(1-\a)B]\timu}\\
\ns\ds={\a[\timu^\top C\timu]C\timu\timu^\top(B-C)\over\timu^\top[\a C+(1-\a)B]\timu}.\ea$$
Hence,
$$\ba{ll}
\ns\ds\h G=\a B+(1-\a)C+\a{(1-\a)B\timu\timu^\top(C-B)\over\timu^\top[\a C+(1-\a)B]\timu}+(1-\a){\a C\timu\timu^\top(B-C)\over\timu^\top[\a C+(1-\a)B]\timu}\\
\ns\ds\q=\a B+(1-\a)C-\a(1-\a){(B-C)\timu\timu^\top(C-B)\over\timu^\top[\a C+(1-\a)B]\timu}.
\ea$$
This means that for any $x\in E_k$,
\bel{E415}\ba{ll}
\ns\ds G(x)=\a B(x)+(1-\a)C(x)-\a(1-\a){\big(C(x)-B(x)\big)\m_k
\m_k^\top\big(C(x)-B(x)\big)\over\m_k^\top\big(\a C(x)+(1-\a)B(x)\big)\m_k}\\
\ns\ds\qq\;=\a B(x)+(1-\a)C(x)-\a(1-\a){\big(C(x)-B(x)\big)\ell(x)\ell(x)^\top\big(C(x)-B(x)\big)\over\ell(x)^\top
\big(\a C(x)+(1-\a)B(x)\big)\ell(x)}.\ea\ee
The proof is completed. \endpf

\vskip-0.5cm

By Lemma \ref{T401}, we know that the function $G:\O\to\dbS^n_+$ appears in the above lemma satisfies the following:
$$\l|\xi|^2\les\lan G(x)\xi,\xi\ran\les\L|\xi|^2,\qq\forall x\in\O,~\xi\in S^{n-1}.$$

\bl{T404} \sl Let {\rm(S1)--(S5)} hold and $(\bar y(\cd),\bar u(\cd))$ be an optimal pair of Problem {\rm(C)}. Let $u(\cd)\in\cU$, $\m_1,\m_2,\ldots,\m_m\in \dbR^n\setminus\{0\}$ and
$E_1,E_2,\cds,E_m$ be mutually disjoint measurable sets such that $\ds\bigcup^m_{i=1}E_i=\O$. Define
$$\ell(\cd)=\sum_{k=1}^m\m_k\chi_{E_k}(\cd).$$
Let $\a\in (0,1)$ and $y^\a(\cd\,;\ell(\cd))$ be the weak solution of
\bel{E419}\left\{\2n\ba{ll}
\ns\ds\ds-\Div\big(\h G(x)\nabla y^\a(x;\ell(\cd))\big)=\a f(x,y^\a(x;\ell(\cd)),u(x))+(1-\a)f(x,y^\a(x;\ell(\cd)),u(x)),\qq\eqin \O,\\
\ns\ds y^\a(x;\ell(\cd))=0,\qq x\in\pa\O, \ea\right.\ee
with $\h G(\cd)\in L^\infty(\O;\dbS^n_+)$ being given by
\bel{E420}\ba{ll}
\ns\ds\h G(x)=\a A(x,u(x))+(1-\a)A(x,\bu(x))\\
\ns\ds\qq\q-{\a(1-\a)\big[A(x,u(x))-A(x,\bu(x))\big]\ell(x)\ell(x)^\top
\big[A(x,u(x))-A(x,\bu(x))\big]\over(1-\a)\ell(x)^\top A(x,u(x))\ell(x)+\a\ell(x)^\top
A(x,\bu(x))\ell(x)},\qq\forall x\in\O.\ea\ee
Then
\bel{E421}J^\a(u(\cd),\ell(\cd))\equiv\int_\O\Big(\a
f^0(x,y^\a(x;\ell(\cd)),u(x))+(1-\a)f^0(x,y^\a(x;\ell(\cd)),u(x))\Big)\,dx\ges J(\bar u(\cd)).\ee
\el

\it Proof. \rm We split the proof into two steps.

\ms

\textbf{Step I.} First, let all the components of $\m_1,\m_2,\ldots,\m_m\in\dbQ^n\setminus\{0\}$.
Let $P$ be a positive integer such that all the components of $P\m_1,P\m_2,\ldots,
P\m_m$ are integers. For $\e>0$, define
\bel{E422}u^{\a,\e}(x)=\left\{\2n\ba{ll}
\ns\ds u(x),\qq\hb{if }\big\{\lan{x\over\e},\m_k\ran\big\}\in [0,\a),~x\in E_k,\\
\ns\ds\bar u(x),\qq\hb{if }\big\{\lan{x\over\e},\m_k\ran\big\}\in [\a,1),~x\in E_k.\ea\right.\ee
Let $y^{\a,\e}(\cd)=y(\cd;u^{\a,\e}(\cd))$ be the solution to the state equation \rf{state}
corresponding to the control $u^{\a,\e}(\cd)$. It is standard that $y^{\a,\e}(\cd)$ is
uniformly bounded in $H^1_0(\O)$ and $L^\infty(\O)$. Thus, along a subsequence $\e\da0$,
$y^{\a,\e}(\cd)$ converges to some $z^\a(\cd)$ weakly in $H^1_0(\O)$ and strongly in $L^2(\O)$.
Let $\wt y^{\a,\e}(\cd)\in H^1_0(\O)$ be the weak solution of
\bel{E423}\left\{\2n\ba{ll}
\ns\ds-\Div\big(A(x,u^{\a,\e}(x))\nabla
\wt y^{\a,\e}(x)\big)=\a f(x,z^\a(x),u(x))+(1-\a)f(x,z^\a(x),u(x)),\qq\hb{in }\O,\\
\ns\ds\wt y^{\a,\e}(x)=0,\qq\hb{ on }\pa\O,\ea\right.\ee
Then $\hy^{\a,\e}(\cd)\equiv y^{\a,\e}(\cd)-\wt y^{\a,\e}(\cd)$ satisfies
\bel{E424}\left\{\2n\ba{ll}
\ns\ds-\Div\big(A(x,u^{\a,\e}(x))\nabla\hy^{\a,\e}(x)\big)=f(x,y^{\a,\e}(x),u^{\a,\e}(x))\\
\ns\ds\qq\qq\qq\qq\qq\qq\qq-\a f(x,z^\a(x),u(x))-(1-\a)f(x,z^\a(x),u(x)),\qq\eqin \O,\\
\ns\ds\hy^{\a,\ve}(x)=0,  \qq\eqon \pa\O,\ea\right.\ee
By (S3) and the boundedness of $y^{\a,\ve}(\cd)$ in $L^\infty(\O)$,
we can see that
$$f(\cd\,,y^{\a,\e}(\cd),u^{\a,\e}(\cd))-f(\cd,z^\a(\cd),u^{\a,\e}(\cd))$$
converges strongly in $L^2(\O)$. In addition, by Riemann-Lebesgue's Theorem, we see that
$$f\big(\cd\,,y^\a(\cd\,;\ell(\cd)),u^{\a,\e}(\cd)\big)\to\a f(\cd\,,z^\a(\cd),u(\cd))
+(1-\a)f(\cd\,,z^\a(\cd),u(\cd)),\qq\hb{weakly in $L^2(\O)$ as $\e\to0$.}$$
Then, (S2) and \eqref{E424} imply that $\hy^{\a,\e}(\cd)$ converges weakly to zero in $H^1_0(\O)$ (as $\e\da0$). Consequently, along a subsequence $\e\da0$, $\tiy^{\a,\e}(\cd)$ converges to $z^\a(\cd)$ weakly in $H^1_0(\O)$. Combining this with Lemma \ref{T403}, we get $z^\a(\cd)=y^\a(\cd\,;\ell(\cd))$, where $y^\a(\cd\,;\ell(\cd))$ is the weak solution of \eqref{E419}. This implies $y^{\a,\e}(\cd)$ itself converges to $y^\a(\cd\,;\ell(\cd))$ weakly in $H^1_0(\O)$, strongly in $L^2(\O)$ (as $\e\da0$). Hence,
\bel{E425}\ba{ll}
\ns\ds J(\bu(\cd))\les\lim_{\e\da0}\int_\O f^0\big(x,y^{\a,\e}(x),u^{\a,\e}(x)\big)\,dx=\lim_{\e\da0}\int_\O f^0\big(x,y^\a(x;\ell(\cd)),u^{\a,\e}(x)\big)\,dx\\
\ns\ds\qq\q=\int_\O\Big(\a f^0\big(x,y^\a(x;\ell(\cd)),u(x)\big)+(1-\a)f^0\big(x,y^\a(x;\ell(\cd)),\bar u(x)\big)\Big)\,dx
=J^\a(u(\cd),\ell(\cd)),\ea\ee
proving \eqref{E421} for the case of $\ell(\cd)$ valued in rational numbers.

\ms

\textbf{Step II.} Now, let $\m_1,\cds,\m_m\in\dbR^n\setminus\{0\}$. We can select $\m_1^\l,\m_2^\l,\ldots,\m_m^\l\in\dbQ^n\setminus\{0\}$ such that $\m_k^\l\to\m_k$ as $\l\to+\infty$. By
Step I, one has
\bel{E426}\int_\O\Big(\a f^0(x,y^\a_\l(x),u(x))+(1-\a)f^0(x,y^\a_\l(x),u(x))\Big)\,
dx\ges J(\bar u(\cd)),\ee
where
\bel{E427}\left\{\2n\ba{ll}
\ns\ds-\Div\big(G_\l(x)\nabla y^\a_\l(x)\big)=\a f(x,y^\a_\l(x),u(x))+(1-\a)f(x,y^\a_\l(x),u(x)),\qq\eqin \O,\\
\ns\ds y^\a_\l\big|_{\,\pa\O}=0 \ea\right.\ee
with $G_\l(\cd)\in L^\infty(\O;\dbS^n_+)$ given by
\bel{E428}\ba{ll}
\ns\ds G_\l(x)=\a A(x,u(x))+(1-\a)A(x,\bu(x))\\
\ns\ds\qq\q-{\a(1-\a)\big[A(x,u(x))-A(x,\bu(x))\big]\m^\l_k
(\m^\l_k)^\top\big[A(x,u(x))-A(x,\bu(x))\big] \over (1-\a)
(\m^\l_k)^\top A(x,u(x))\m^\l_k+\a(\m^\l_k)^\top A(x,\bu(x))\m^\l_k},\qq\forall x\in E_k.\ea\ee
By Lemma \ref{T401}, $G_\l(x)\in\dbS^n_+$ for almost all $x\in\O$. Clearly, $G_\l(\cd)$ converges to $G(\cd)$ strongly
in $L^\infty(\O)$. Then by Lemma \ref{T302} and a standard argument, we have the convergence of $y^\a_\l(\cd)$ to $y^\a(\cd)$ strongly in $H^1_0(\O)$. Then \eqref{E421} follows. \endpf

\vskip-0.5cm

We now turn to a proof of Proposition \ref{T303}.

\ms

\it Proof of Proposition \ref{T303}. \rm By Luzin's Theorem, for any integer $\l\ges1$, there exists a closed subset
$F_\l$ of $\bar\O$, such that $\ell(\cd)$ is continuous on $F_\l$ and $|\bar\O\setminus
F_\l|\les{1\over\l}$, where $|S|$ stands for the Lebesgue measure of the set $S$. Since $F_\l$ is also bounded, $\ell(\cd)$ is uniformly continuous on $F_\l$. Thus, there exist disjoint measurable sets $E_{\l 1}$, $E_{\l 2},\ldots,E_{\l m_\l}$
such that
$$\ds\bigcup^{m_\l}_{k=1}E_{\l k}=F_\l,\qq\sup_{x,\tix\in E_{\l k}}|\ell(x)-\ell(\tix)|\les{1\over\l},\qq\forall j=1,2,\ldots,m_\l.$$
Choosing arbitrary $x^{\l k}$ from $E_{\l k}$ and $x^{\l0}$ from $E_{\l0}\equiv\bar\O\setminus F_\l$, we define
$$\m_{\l k}=\ell(x^{\l k}),\qq k=0,1,2,\ldots,m_\l.$$
By Lemma \ref{T403},
\bel{E429}\int_\O\Big(\a f^0(x,y^\a_\l(x),u(x))+(1-\a)f^0(x,y^\a_\l(x),u(x))\Big)\,dx\ges
J(\bar u(\cd)),\ee
where $y^\a_\l(\cd)$ is the weak solution of the following:
\bel{E430}\left\{\2n\ba{ll}
\ns\ds-\Div\big(G^\a_\l(x)\nabla y^\a_\l(x)\big)=\a f(x,y^\a_\l(x),u(x))
+(1-\a)f(x,y^\a_\l(x),u(x)),\qq\eqin \O,\\
\ns\ds y^\a_\l\big|_{\,\pa\O}=0,\ea\right.\ee
with $G^\a_\l(\cd)\in L^\infty(\O;\dbS^n_+)$ being given by
\bel{E431}\ba{ll}
\ns\ds G^\a_\l(x)=\a A(x,u(x))+(1-\a)A(x,\bar u(x))\\ [2mm]
\ns\ds\qq\qq\q-{\a(1-\a)\big[A(x,u(x))-A(x,\bar u(x))\big]\m_{\l k}
\m_{\l k}^\top\big[A(x,u(x))-A(x,\bar u(x))\big]\over(1-\a)\m_{\l k}^\top
A(x,u(x))\m_{\l k}+\a\m_{\l k}^\top A(x,\bar u(x))\m_{\l k}},\\ [3mm]
\ns\ds\qq\qq\qq\qq\qq\qq\qq\forall x\in E_{\l k},~k=0,1,2,\ldots,m_\l.\ea\ee
Obviously, $G^\a_\l(\cd)$ converges to $A^\a(\cd)$ strongly in $L^2(\O)$, as $\l\to\infty$,
where $A^\a(\cd)$ is defined  by \rf{A^a(ell)}. Thus it follows from (S2)--(S3) and Lemma
\ref{T302} that $y^\a_\l(\cd)$ converges to $y^\a(\cd\,;\ell(\cd))$, as $\l\to\infty$, strongly in $H^1_0(\O)$.
We then obtain \eqref{J^a(u,ell)-J>0}. \endpf

\section{Second-Order Necessary Conditions.}

In this section, we are going to prove Theorem \ref{T204}. For readers' convenience, we will rewrite the relevant equations when needed. We first establish the following lemma.

\bl{T501} \sl Let {\rm(S1)--(S5)} hold. Let $\bu(\cd),u(\cd)\in\cU$ and $\ell(\cd)\in\cL$. For any $\a\in(0,1)$, define
\bel{A^a(ell)*}\ba{ll}
\ns\ds A^\a(x;\ell(\cd))=\a A(x,u(x))+(1-\a)A(x,\bu(x))\\
\ns\ds\qq\qq\qq-{\a(1-\a)\big[A(x,u(x))-A(x,\bu(x))\big]\ell(x) \ell(x)^\top
\big[A(x,u(x))-A(x,\bu(x))\big] \over (1-\a) \ell(x)^\top
A(x,u(x))\ell(x)+\a \ell(x)^\top A (x,\bu(x))\ell(x)}.\ea\ee
Let $\by(\cd)$, $y^\a(\cd\,;\ell(\cd))$ be the weak solutions of the following equations:
\bel{state*}\left\{\2n\ba{ll}
\ns\ds-\Div\big(A(x,\bu(x))\nabla
\by(x)\big)=f(x,\by(x),\bu(x)),\qq\hb{in }~\O,\\
\ns\ds\by(x)=0,\qq\qq\qq\qq\qq\hb{on }~\pa\O,\ea\right.\ee
and
\bel{state-y^a(ell)*}\left\{\2n\ba{ll}
\ns\ds-\Div\(A^\a(x;\ell(\cd))\nabla y^\a(x;\ell(\cd))\)=(1-\a)f\big(x,y^\a(x;\ell(\cd)),\bar u(x)\big)\1n+\1n\a f\big(x,y^\a(x;\ell(\cd)),u(x)\big),\q\eqin\O,\\
\ns\ds y^\a(x;\ell(\cd))=0,\qq\qq\eqin\pa\O,\ea\right.\ee
respectively. Then, as $\a\da0$, $\ds Y^\a\equiv {y^\a(\cd\,;\ell(\cd))-\by(\cd)\over \a}$ converges to $Y(\cd)$ weakly in $H^1_0(\O)$ with $Y(\cd)$ being the weak solution of
\bel{variation*}\left\{\2n\ba{ll}
\ns\ds-\Div\big(A(x,\bu(x))\nabla Y(x)\big)=f_y(x,\by(x),\bu(x))\,Y(x)+\Div \big(\Th(x)\nabla\by(x)\big)\\
\ns\ds\qq\qq\qq\qq\qq\qq\q+f(x,\by(x),u(x))-f(x,\by(x),\bu(x)),\qq\eqin  \O,\\
\ns\ds Y\big|_{\pa\O}=0,\ea\right.\ee
where
\bel{Th*}\ba{ll}
\ns\ds\Th(x)=A(x,u(x))-A(x,\bu(x))\\
\ns\ds\qq\q-{\big[A(x,u(x))-A(x,\bu(x))\big]\ell(x)\ell(x)^\top\big[A(x,u(x))-A(x,\bu(x))\big]\over
\ell(x)^\top A(x,u(x))\ell(x)}\,.\ea\ee
\el

\it Proof. \rm Let $y^\a(\cd)=y^\a(\cd\,;\ell(\cd))$. We have
\bel{E506}\ba{ll}
\ns\ds-\Div\big(A(x,\bu(x))\nabla Y^\a(x)\big)\\
\ns\ds={f(x,y^\a(x),\bu(x))-f(x,\by(x),\bu(x))\over \a}+
f(x,y^\a(x),u(x))-f(x,y^\a(x),\bu(x))\\
\ns\ds\qq-\Div\Big({A(x,\bu(x))-A^a(x)\over \a}\nabla y^\a(x)\Big)\\
\ns\ds=\(\int^1_0  f_y(x,\by(x)+t(y^\a(x)-\by(x)),\bu(x))\, dt\)Y^\a(x)\\
\ns\ds\qq+f(x,y^\a(x),u(x))-f(x,y^\a(x),\bu(x)) +\Div\Big((A(x,u(x))-A(x,\bu(x)))\nabla y^\a(x)\Big)\\
\ns\ds\qq-(1-\a)\Div\Big({ \big[A(x,u(x))-A(x,\bu(x))\big]\ell(x)
\ell(x)^\top \big[A(x,u(x))-A(x,\bu(x))\big] \over (1-\a)
\ell(x)^\top A(x,u(x))\ell(x)+\a
\ell(x)^\top A (x,\bu(x))\ell(x)}\nabla y^\a(x)\Big).\ea\ee
By Proposition \ref{T303} and (S2)--(S3), one can see that $y^\a(\cd)$ is bounded uniformly in $H^1_0(\O)$. Thus, we can prove that $Y^\a(\cd)$ is bounded uniformly in $H^1_0(\O)$, and as $\a\da0$, $y^\a(\cd)$ converges to $\by(\cd)$ strongly in $H^1_0(\O)$, $Y^\a(\cd)$ converges to $Y(\cd)$ weakly in $H^1_0(\O)$ with $Y(\cd)$ being the weak solution of \eqref{variation*}. \endpf

\vskip-0.5cm

Let us define
\bel{Z^a}Z^\a(\cd)\equiv {Y^\a(\cd)-Y(\cd)\over\a}.\ee
Then it is natural to expect that as $\a\da0$, $Z^\a(\cd)\to Z(\cd)$, weakly in $H^1_0(\O)$, with $Z(\cd)$ being the weak solution of the following equation:
\bel{E508}\left\{\2n\ba{ll}
\ns\ds-\Div\big(A(x,\bu(x))\nabla Z(x)\big)=f_y(x,\by(x),\bu(x))Z(x)+{1\over2}f_{yy}(x,\by(x),\bu(x))|Y (x)|^2\\
\ns\ds\qq\qq\qq\qq\qq\q+\big(f_y(x,\by(x),u(x))-f_y(x,\by(x),\bu(x))\big)Y(x)\\
\ns\ds\qq\qq\qq\qq\qq\q+\Div\Big((A(x,u(x))-A(x,\bu(x)))\nabla Y (x)\Big)+\Div\Big(\Upsilon(x)\nabla \by (x)\Big),\qq\hb{in  }\O,\\
\ns\ds Z\big|_{\pa\O}=0,\ea\right.\ee
where
$$\Upsilon(x)=\big[A(x,u(x))-A(x,\bu(x))\big]\ell(x) \ell(x)^\top
\big[A(x,u(x))-A(x,\bu(x))\big]{\ell(x)^\top
A(x,\bu(x))\ell(x)\over  [\ell(x)^\top A(x,u(x))\ell(x)]^2}.$$
However, this seems not to be true when $n$ is a large. The main reason is that  $A(\cd,\bu(\cd))$ is only bounded and measurable so that we have $W^{1,p}(\O)$ estimates of $\ds Z^\a(\cd)$ only for $p$ near $2$. Thus, we do not have a weak maximum principle for $Y(\cd)$, i.e., we usually cannot
get uniform boundedness for $Y(\cd)$ unless $n\les2$. Thus, generally, $Y(\cd)^2$ times an $L^\infty$ function might not be in $H^{-1}(\O)$ unless $n\les6$. Hence, $Z(\cd)$ is probably not well-defined for large $n$. Fortunately, we have the following relation, which will be sufficient for the proof of our main result---Theorem \ref{T204}. Recalling the definition of $Z^\a(\cd)$, we have
\bel{E507}\ba{ll}
\ns\ds-\Div\big(A(\bu)\nabla Z^\a\big)\\
\ns\ds={1\over\a}\[\(\int^1_0f_y\big(\by+t(y^\a-\by),\bu\big)\,dt\)Y^\a
-f_y(\by,\bu)\,Y\]+{f(y^\a,u)-f(\by,u)\over\a}-{f(y^\a,\bu)-f(\by,\bu)\over\a}\\
\ns\ds\q+{1\over\a}\[\Div\((A(u)-A(\bu)\nabla y^\a\)-(1-\a)\Div\({\big[A(u)-A(\bu)
\big]\ell\ell^\top\big[A(u)-A(\bu)\big]\over(1-\a)\ell^\top A(u)\ell+\a
\ell^\top A(\bu)\ell}\nabla y^\a\)-\Div\big(\Th\nabla\by\big)\]\\
\ns\ds=f_y(\by,\bu)\,Z^\a+\(\int^1_0{f_y(\by+t(y^\a-\by),\bu)
-f_y(\by,\bu)\over\a}\,dt\)Y^\a\\
\ns\ds\q+{f(y^\a,u)-f(\by,u)\over\a}-{f(y^\a,\bu)-f(\by,\bu)\over\a}+\Div\((A(u)-A(\bu))
\nabla Y^\a\)\\
\ns\ds\qq+\Div\({\big[A(u)-A(\bu)\big]\ell\ell^\top\big[A(u)-A(\bu)\big]\over
(1-\a)\ell^\top A(u)\ell+\a\ell^\top A(\bu)\ell}\nabla y^\a\)\\
\ns\ds\qq-\Div\Big({\ell^\top\big[A(u)-A(\bu)\big]\ell\over\ell^\top A(u)\ell[(1-\a)
\ell^\top A(u)\ell+\a\ell^\top A(\bu)\ell]}\big[A(u)-A(\bu)\big]\ell\ell^\top
\big[A(u)-A(\bu)\big]\nabla y^\a\Big)\\
\ns\ds=f_y(\by,\bu)\,Z^\a+\(\int^1_0 dt\int^1_0 tf_{yy}(\by+\t t(y^\a-\by),\bu)\,d\t\)|Y^\a|^2\\
\ns\ds\qq+\[\int^1_0\(f_y(\by+t(y^\a-\by),u)-f_y(\by+t(y^\a-\by),\bu)\)dt\]Y^\a+\Div\Big((A(u)-A(\bu)\nabla Y^\a\Big)\\
\ns\ds\qq+\Div\Big({\ell^\top A(\bu)\ell\over\ell^\top A(u)\ell
[(1-\a)\ell^\top A(u)\ell+\a\ell^\top A(\bu)\ell]}\big[A(u)-A(\bu)\big]\ell\ell^\top
\big[A(u)-A(\bu)\big]\nabla y^\a\Big).\ea\ee

Now, we are ready to prove our main result.

\ms

\it Proof of Theorem \ref{T204}. \rm Let $(u(\cd),\ell(\cd))\in\cV_0(\bar u(\cd))$ and $\a\in
(0,1)$. Let $y^\a(\cd)$ be the weak solution to \rf{state-y^a(ell)}, with $A^\a(\cd)$ and $J^\a(u(\cd),\ell(\cd))$ be defined by \rf{A^a(ell)} and \rf{J^a(u,ell)}, respectively. Let $\ds Y^\a(\cd)={y^\a(\cd)-\by(\cd)\over
\a}$. Then $\ds Y^\a(\cd)$ satisfies \eqref{E506}. We have
\bel{E509}\ba{ll}
\ns\ds\int_\O\(f^0(\by,u)-f^0(\by,\bu)+f^0_y(\by,\bu)Y\)dx\\
\ns\ds=\int_\O\(f^0(\by,u)-f^0(\by,\bu)\)\,dx+\int_\O\(f_y(\by,\bu)\,\bar\psi+
\Div\big(A(\bu)\nabla\bar\psi\big)Y\)dx\\
\ns\ds=\int_\O\(f^0(\by,u)-f^0(\by,\bu)+f_y(\by,\bu)Y+\Div\big(A(\bu)\nabla Y\big)\)
\bar\psi\, dx\\
\ns\ds=\int_\O\[f^0(\by,u)-f^0(\by,\bu)-\(f(\by,u)-f(\by,\bu)+\Div\big(
\Th\nabla\by\big)\)\bar\psi\]dx\\
\ns\ds=\2n\int_\O\2n\Big\{H\big(\by,\bar\psi,\nabla\by,\nabla\bar\psi,\bu\big)\1n
-\1n H\big(\by,\bar\psi,\nabla\by,\nabla\bar\psi,u\big)\1n-\1n{\lan[A(u)\1n-\1n A(\bu)]\nabla\bar y,\ell\ran\,\lan\big[A(u)\1n-\1n A(\bu)\big]\nabla\bar\psi,\ell\ran\over\lan A(u)\ell,\ell\ran}\Big\}dx\1n=\1n0,\ea\ee
where $\Th(\cd)$ is given by \eqref{Th*}. In fact,
$$\int_\O\(f^0(\by,u)-f^0(\by,\bu)+f^0_y(\by,\bu)Y(x)\)dx=\lim_{\a\da0}{J^\a(u(\cd),
\ell(\cd))-J(\bu(\cd))\over\a}=0,$$
due to the partial singularity of $\bar u(\cd)$ at $(u(\cd),\ell(\cd))\in\cV_0(\bar u(\cd))$.
This leads to
\bel{}\int_\O\(f^0(\bar y,u)-f^0(\bar y,\bar u)\)dx=-\int_\O f^0_y(\bar y,\bar u)Y dx.\ee

Then, using \eqref{E509}, \eqref{adjoint} and \eqref{E506}, we have (suppressing $x$ whenever no confusion would be caused, for notational simplicity)
\bel{E510}\ba{ll}
\ns\ds J^\a(u(\cd),\ell(\cd))-J(\bu(\cd))\equiv\a\int_\O f^0(y^\a,u)dx+(1-\a)\int_\O f^0(y^\a,\bar u)dx-\int_\O f^0(\bar y,\bar u)dx \\
\ns\ds=\a\int_\O\(f^0(y^\a,u)-f^0(y^\a,\bu)\)\,dx+\int_\O\(f^0(y^\a,\bu)-f^0(\by,\bu)\)\,dx \\
\ns\ds=\a\int_\O\(f^0(y^\a,u)-f^0(\bar y,u)-f^0(y^\a,\bar u)+f^0(\bar y,\bar u)+f^0(\bar y,u)-f^0(\bar y,\bar u)\)\,dx\\
\ns\ds\q+\int_\O\(f^0(y^\a,\bar u)-f^0(\bar y,\bar u)\)dx\\
\ns\ds=\a^2\int_\O\[\int^1_0\(f^0_y(\by+\a t Y^\a,u)-f^0_y(\by+\a t Y^\a,\bu)\)dt\]Y^\a\,dx\\
\ns\ds\q+\a\int_\O\[\(\int_0^1f^0_y(\bar y+\a tY^\a,\bar u)dt\)Y^\a-f^0_y(\bar y,\bar u)Y\]dx\\
\ns\ds=\a^2\int_\O\[\int^1_0\(f^0_y(\by+\a t Y^\a,u)-f^0_y(\by+\a t Y^\a,\bu)\)dt\]Y^\a\,dx\\
\ns\ds\q+\a\int_\O\[\(\int_0^1f^0_y(\bar y+\a tY^\a,\bar u)dt-f^0_y(\bar y,\bar u)\)Y^\a
+\a f^0_y(\bar y,\bar u)Z^\a\]dx\\
\ns\ds=\a^2\int_\O\[\int^1_0\(f^0_y(\by+\a t Y^\a,u)-f^0_y(\by+\a t Y^\a,\bu)\)dt\]Y^\a\,dx\\
\ns\ds\q+\a^2\2n\int_\O\1n\[\1n\int^1_0\2n\(\1n\int^1_0 tf^0_{yy}(\by+\a t\t Y^\a,\bu)d\t\)dt\]|Y^\a|^2dx\1n+\1n\a^2\2n\int_\O\1n\(f_y(\by,\bu)\,\bar\psi\1n+\2n\Div\big(A(\bu)\nabla\bar\psi\big)\)Z^\a\,dx \\
\ns\ds=\a^2\int_\O\[\int^1_0\(f^0_y(\by+\a t Y^\a,u)-f^0_y(\by+\a t Y^\a,\bu)\)dt\]Y^\a\,dx\\
\ns\ds\q+\a^2\2n\int_\O\1n\[\1n\int^1_0\2n\(\1n\int^1_0 tf^0_{yy}(\by+\a t\t Y^\a,\bu)d\t\)dt\]|Y^\a|^2dx\1n+\1n\a^2\2n\int_\O\1n\(f_y(\by,\bu)Z^\a\1n+\2n\Div\big(A(\bu)
\nabla Z^\a\big)\)\bar\psi\,dx \\
\ns\ds=\a^2\Big\{\int_\O\[\int^1_0\(f^0_y(\by+\a t Y^\a,u)-f^0_y(\by+\a t Y^\a,\bu)\)dt\]Y^\a\,dx\\
\ns\ds\q+\2n\int_\O\1n\[\1n\int^1_0\2n\(\1n\int^1_0 tf^0_{yy}(\by+\a t\t Y^\a,\bu)d\t\)dt\]|Y^\a|^2dx\\
\ns\ds\q-\int_\O\[\int^1_0\(\int^1_0tf_{yy}(\by+\a\t tY^\a,\bu)d\t\)dt\]\bar\psi|Y^\a|^2dx\\
\ns\ds\q-\int_\O\[\int^1_0\(f_y(\by+\a tY^\a,u)-f_y(\by+\a tY^\a,\bu)\)dt\]\bar\psi\,Y^\a dx+\int_\O\lan(A(u)-A(\bu))\nabla Y^\a,\nabla\bar\psi\ran\,dx\\
\ns\ds\q+\int_\O\Big\langle{\ell^\top A(\bu)\ell\big[A(u)-A(\bu)\big]\ell\ell^\top
\big[A(u)-A(\bu)\big]\nabla y^\a\over\ell^\top A(u)\ell[(1-\a)\ell^\top A(u)\ell+\a\ell^\top A(\bu)\ell]},\nabla\bar\psi\Big\rangle\,dx\Big\}\\
\ns\ds=\a^2\Big\{\int_\O\[\int^1_0\(
H_y(\by+\a tY^\a,\bar\psi,\nabla\by,\nabla\bar\psi,\bu)-H_y(\by+\a tY^\a,\bar\psi,\nabla \by,\nabla \bar\psi,u)\)dt\]Y^\a dx\\
\ns\ds\q-\1n\int_\O\[\int^1_0\2n\(\int^1_0tH_{yy}(\by+\a t\t Y^\a,\bar\psi,\nabla\by,\nabla\bar\psi,\bu)d\t\)dt\]|Y^\a|^2dx\1n+\2n\int_\O\2n\lan(A(u)-A(\bu))\nabla Y^\a,\nabla\bar\psi\ran\,dx\\
\ns\ds\q+\int_\O\Big\langle{[\ell^\top A(\bu)\ell]\{\ell^\top\big[A(u)-A(\bu)\big]\nabla\bar\psi\}\{\ell^\top
\big[A(u)-A(\bu)\big]\nabla y^\a\}\over\ell^\top A(u)\ell[(1-\a)\ell^\top A(u)\ell+\a\ell^\top A(\bu)\ell]}dx\Big\}.\ea\ee
Thus,
$$\ba{ll}
\ns\ds0\les\lim_{\a\da0}{J^\a(u(\cd),\ell(\cd))-J(\bu(\cd))\over\a^2} \\
\ns\ds=\int_\O\[\Big(H_y(x,\by(x),\bar\psi(x),\nabla\by(x),\nabla\bar\psi(x),\bu(x))- H_y(x,\by(x),\bar\psi(x),\nabla\by(x),\nabla\bar\psi(x),u(x))\Big)Y(x)\\
\ns\ds\q-{1\over2}H_{yy}(x,\by(x),\bar\psi(x),\nabla\by(x),\nabla\bar\psi(x),\bu(x))|Y(x)|^2
+\lan\big(A(x,u(x))-A(x,\bu(x)\big)\nabla\bar\psi(x),\nabla Y(x)\ran \\
\ns\ds\q+\ell(x)^\top \big[A(x,u(x))-A(x,\bu(x))\big]\nabla\bar\psi(x)\ell(x)^\top \big[A(x,u(x))-A(x,\bu(x))\big]\nabla\by(x){\big[\ell(x)^\top A(x,\bu(x))\ell(x)\big]\over
[\ell(x)^\top A(x,u(x))\ell(x) ]^2}\Big\}dx\\
\ns\ds=\int_\O\[\Big(H_y(x,\by(x),\bar\psi(x),\nabla\by(x),\nabla\bar\psi(x),\bu(x))- H_y(x,\by(x),\bar\psi(x),\nabla\by(x),\nabla\bar\psi(x),u(x))\Big)Y(x)\\
\ns\ds\q-{1\over2}H_{yy}(x,\by(x),\bar\psi(x),\nabla\by(x),\nabla\bar\psi(x),\bu(x))|Y(x)|^2
+\lan\big(A(x,u(x))-A(x,\bu(x)\big)\nabla\bar\psi(x),\nabla Y(x)\ran \\
\ns\ds\q+\1n\(\1n H(x,\by(x),\bar\psi(x),\nabla\by(x),\nabla\bar\psi(x),\bu(x))\1n-\1n H(x,\by(x),\bar\psi(x),\nabla\by(x),\nabla
\bar\psi(x),u(x))\1n\){\ell(x)^\top A(x,\bu(x))\ell(x)\over
 \ell(x)^\top A(x,u(x))\ell(x)}\Big\}dx.\ea$$
This completes the proof. \endpf

\section{Concluding Remarks.}

We have established the second-order necessary conditions for the optimal controls of Problem (C). There are some challenging problems left open. We list some of them here, for which we are still working on with our great efforts.

\ms

$\bullet$ Construction of suitable examples for which our second-necessary conditions could lead to some optimal solutions.

\ms

$\bullet$ The second-order necessary conditions that we obtained looks complicated. Is it possible to have some better forms?

\ms

$\bullet$ Extension to fully non-linear equations.

\ms

We hope to be able to report some further results before long. Also, any participation of other interested researchers are welcome.

\end{document}